\documentclass[11pt,reqno]{amsart}
\usepackage{amsmath}
\usepackage{amssymb, amsthm, amsfonts, mathrsfs}
\usepackage{bm}

\newtheorem*{notation}{Notation}

\newcommand{\al}{\alpha}
\newcommand{\be}{\beta}
\newcommand{\ga}{\gamma}
\newcommand{\de}{\delta}

\newcommand{\rig}{\rightarrow}
\newcommand{\mrig}{\mathrel{-\!\!\!\!\!\rightarrow}}

\newcommand{\Rig}{\Rightarrow}

\newcommand{\bcdw}{\mathbin{\boldsymbol\cdot}}
\newcommand{\bcdn}{\mbox{\boldmath{$\cdot$}}}

\newcommand{\seteq}{\mathrel{\mbox{\,:\!}=\nolinebreak }\,}

\newcommand{\up}{\mathop{\,\uparrow}}
\newcommand{\down}{\mathop{\,\downarrow}}

\newcommand{\sbA}{{\boldsymbol{A}}}
\newcommand{\sbB}{{\boldsymbol{B}}}
\newcommand{\sbC}{{\boldsymbol{C}}}
\newcommand{\sbD}{{\boldsymbol{D}}}
\newcommand{\sbE}{{\boldsymbol{E}}}
\newcommand{\sbF}{{\boldsymbol{F}}}
\newcommand{\sbG}{{\boldsymbol{G}}}

\newcommand{\sbK}{{\boldsymbol{K}}}
\newcommand{\sbL}{{\boldsymbol{L}}}

\newcommand{\sbP}{{\boldsymbol{P}}}

\newcommand{\sbS}{{\boldsymbol{S}}}
\newcommand{\sbT}{{\boldsymbol{T}}}
\newcommand{\sbU}{{\boldsymbol{U}}}
\newcommand{\sbV}{{\boldsymbol{V}}}
\newcommand{\sbW}{{\boldsymbol{W}}}
\newcommand{\sbX}{{\boldsymbol{X}}}
\newcommand{\sbY}{{\boldsymbol{Y}}}
\newcommand{\sbZ}{{\boldsymbol{Z}}}

\newcommand{\tcg}{{\Theta}}

\newcommand{\ov}{\overline}

\let\class=\mathsf                              
\bmdefine{\Fm}{Fm}                                
\bmdefine{\A}{A}                                   
\bmdefine{\B}{B}
\bmdefine{\C}{C}
\bmdefine{\D}{D}

\newcommand{\VVV}{\mathbb{V}}
\newcommand{\QQQ}{\mathbb{Q}}
\newcommand{\OOO}{\mathbb{O}}
\newcommand{\EEE}{\mathbb{E}}

\newcommand{\HHH}{\mathbb{H}}

\newcommand{\PPP}{\mathbb{P}}

\newcommand{\PPU}{\mathbb{P}_{\!\textsc{u}}^{}}
\newcommand{\RRU}{\mathbb{R}_{\textsc{u}}^{}}
\newcommand{\PPS}{\mathbb{P}_{\!\textsc{s}}^{}}
\newcommand{\SSS}{\mathbb{S}}
\newcommand{\III}{\mathbb{I}}
\newcommand{\UUU}{\mathbb{U}}

\bmdefine{\btau}{\tau}                                  
\bmdefine{\brho}{\rho}                                  

\newcommand{\aand}{\mathbin{\,\&\,}}

\theoremstyle{theorem}
\newtheorem{Theorem}{Theorem}[section]
\newtheorem{Lemma}[Theorem]{Lemma}
\newtheorem{Corollary}[Theorem]{Corollary}
\newtheorem{Proposition}[Theorem]{Proposition}
\newtheorem{Fact}[Theorem]{Fact}

\theoremstyle{definition}
\newtheorem{law}[Theorem]{Definition}
\newtheorem{exa}[Theorem]{Example}
\newtheorem{exas}[Theorem]{Examples}
\newtheorem{rema}[Theorem]{Remark}
\newtheorem{Note}[Theorem]{Note}

\theoremstyle{remark}

\begin{document}
\title[Singly Generated Quasivarieties and Residuated Structures]{Singly Generated Quasivarieties\\and Residuated Structures}
\author{T.\ Moraschini}
\address{Institute of Computer Science, Academy of Sciences of the Czech Republic, Pod Vod\'{a}renskou v\v{e}\v{z}\'{i} 2, 182 07 Prague 8, Czech Republic.}
\email{moraschini@cs.cas.cz}
\author{J.G.\ Raftery}
\address{Department of Mathematics and Applied Mathematics,
 University of Pretoria,
 Private Bag X20, Hatfield,
 Pretoria 0028, and DST-NRF Centre of Excellence in Mathematical and Statistical Sciences (CoE-MaSS), South Africa}
\email{{james.raftery@up.ac.za}}
\author{J.J.\ Wannenburg}
\address{Department of Mathematics and Applied Mathematics,
 University of Pretoria,
 Private Bag X20, Hatfield,
 Pretoria 0028, and DST-NRF Centre of Excellence in Mathematical and Statistical Sciences (CoE-MaSS), South Africa}
\email{{jamie.wannenburg@up.ac.za}}
\keywords{Joint embedding property, (passive) structural completeness, retract, De Morgan monoid, Dunn monoid,
Brouwerian algebra.\vspace{1mm}
\vspace{1mm}}
\thanks{This work received funding from the European Union's Horizon 2020 research and innovation programme under the Marie Sklodowska-Curie grant agreement
No~689176 (project ``Syntax Meets Semantics: Methods, Interactions, and Connections in Substructural logics"). The first author was supported by project CZ.02.2.69/0.0/0.0/17\_050/0008361, OPVVV M\v{S}MT, MSCA-IF Lidsk\'{e} zdroje v teoretick\'{e} informatice, and by RVO 67985807.
%
%
The second author was supported in part
by the National Research Foundation of South Africa (UID 85407).
The third author was supported by the DST-NRF Centre of Excellence in Mathematical and Statistical Sciences (CoE-MaSS), South Africa.
The
University of Pretoria Staff Exchange Bursary Programme and
the
CoE-MaSS (RT18ALG/006)
partially funded the first author's travel to Pretoria in 2017 and 2018, respectively.
Opinions expressed and conclusions arrived at are those of the authors and are not necessarily to be attributed to the CoE-MaSS}


\makeatletter
\renewcommand{\labelenumi}{\text{(\theenumi)}}
\renewcommand{\theenumi}{\roman{enumi}}
\renewcommand{\theenumii}{\roman{enumii}}
\renewcommand{\labelenumii}{\text{(\theenumii)}}
\renewcommand{\p@enumii}{\theenumi(\theenumii)}
\makeatother


%


\begin{abstract}
A quasivariety $\class{K}$ of algebras has the joint embedding property (JEP)
iff it is generated by a single algebra $\sbA$.  It is structurally complete iff the free
$\aleph_0$-generated algebra in $\class{K}$ can serve as $\sbA$.  A consequence of this demand,
called `passive structural completeness' (PSC),
is that the nontrivial members of $\class{K}$ all satisfy the same existential positive sentences.
We prove that if $\class{K}$ is PSC then it still has the JEP, and if it has the JEP and its nontrivial
members lack trivial subalgebras, then
its relatively simple members all belong to the universal class generated by one of them.
Under these conditions, if $\class{K}$ is
relatively semisimple then it is generated by one $\class{K}$-simple algebra.  It is a minimal
quasivariety if, moreover, it is PSC but fails to unify some finite set of equations.
We also prove that a quasivariety of finite type, with a finite nontrivial member, is PSC iff its
nontrivial members have a common retract.
The theory is then applied to the
variety of De Morgan monoids, where we isolate the sub(quasi)varieties that are PSC and those that
have the JEP, while throwing fresh light on those that are structurally complete.
The results illuminate the extension lattices of intuitionistic and relevance logics.
\end{abstract}

$\mbox{}$\\{\vspace{-11mm}}

\maketitle

\vspace{-5mm}

\section{Introduction}


Familiar logics often have an algebraic counterpart that is a quasivariety $\class{K}$ of algebras; in
many cases it is a variety.  In this situation, the derivable inference rules of the logic may or may
not be determined by a single set of `truth tables', i.e., by the operation tables of a
single algebra $\sbA\in\class{K}$.  If some member of $\class{K}$ determines the
\emph{finite} rules of the logic, then another
member determines \emph{all} of the
rules (see Remark~\ref{finite to infinite}), so what is needed is only that $\class{K}$ be generated by a single algebra.
Even when $\class{K}$ is a variety, it must be generated \emph{as a quasivariety}
by one of its members, if the generator is to determine rules (as opposed to theorems only).

Obviously, classical propositional logic has this property:
its algebraic counterpart---the variety of Boolean algebras---is generated as a quasivariety
by its unique two-element member.  More surprisingly, the same holds for intuitionistic propositional
logic (though not with a finite algebra), and for the relevance logic $\mathbf{R}$
\cite{Tok79}, but not for its conservative expansion $\mathbf{R^t}$ (with
the so-called Ackermann constants of \cite{AB75}).  In the intuitionistic case, the algebra determining the (possibly
infinite) rules cannot be countable \cite{Wro74b}.

Maltsev \cite{Mal66} proved that a quasivariety $\class{K}$ is generated by a single algebra iff it has the
\emph{joint embedding property} (JEP), i.e., any two nontrivial members of $\class{K}$ can both be embedded
into some third member.  \L os and Suszko \cite{LS58} characterized this demand by a syntactic `relevance
principle' (finitized in Definition~\ref{relevance principle} below).  Various strengthenings of the JEP
have received attention in the literature.  Their names reflect logical origins, but we choose maximally
transparent characterizations here as definitions.  One such strengthening, called \emph{structural completeness},
asks (in effect) that a quasivariety be generated by its free $\aleph_0$-generated member.  A
weaker variant, now called \emph{passive structural completeness} (PSC), amounts to the demand that
any two nontrivial members of $\class{K}$ have the same existential positive theory.  This hereditary
property still implies the JEP
(Theorem~\ref{Thm:PSC-implies-JEP}).

Our original goal was to investigate these properties for classes of De Morgan monoids (i.e., the
models of $\mathbf{R^t}$).  It became clear, however, that in many of our results, large parts of the
proofs had a general universal algebraic (or even model-theoretic) character, so the first half of this
paper concerns such generalities.  We call $\class{K}$ a \emph{Koll\'{a}r quasivariety} (after \cite{Kol79})
if its nontrivial members lack trivial subalgebras.  We prove that if such a quasivariety has the JEP,
then its relatively simple members all belong to the universal class generated by one of them (Theorem~\ref{Thm:simple}).
If, in addition, $\class{K}$ is relatively semisimple, then it is generated (as a quasivariety) by one
$\class{K}$-\emph{simple} algebra.  We prove that a quasivariety of finite type with a finite nontrivial member is PSC iff its nontrivial
members have a common retract (Theorem~\ref{Thm:retract}).

\enlargethispage{5pt}

The second half of the paper deals with (quasi)varieties of De Morgan monoids, and some closely related residuated
structures.  Among other results, we describe completely the varieties of De Morgan monoids
that are PSC (Theorem~\ref{psc main}), and characterize those with the JEP (Theorem~\ref{my conjecture}).  The
structurally complete varieties of De Morgan monoids fall into two classes---a denumerable family that is fully
transparent and a more opaque collection of subvarieties of a certain PSC variety $\class{M}$.  We prove that this $\class{M}$
also has uncountably many structurally incomplete subvarieties, by exhibiting $2^{\aleph_0}$ structurally incomplete
varieties of Brouwerian algebras (of depth $3$) and applying a `reflection' construction (Theorems~\ref{m continuum not sc} and
\ref{continuum brouwerian not sc}).

\section{Preliminaries}

We deal with structures $\sbA=\langle A;F;R\rangle$, where $F$ [resp.\ $R$] is a set of finitary
operations [resp.\ relations] on the non-empty set $A$.  Constants, i.e., distinguished elements
of $\sbA$, are treated as nullary basic operations.  We call $\sbA$ \emph{finite} [resp.\ \emph{trivial}]
if its \emph{universe} $A$ is finite [resp.\ $\left|A\right|=1$ and $R$ consists of non-empty relations].
Of course, structures model first-order signatures, a.k.a.\ types.
An infinite
set $\mathit{Var}$ of variables is fixed for the entire discussion.
Formulas of any kind are assumed to involve only variables from $\mathit{Var}$.
Recall that a first-order formula with no free variable is called a \emph{sentence}.
For a set $\Sigma$ of
first-order formulas,
the notation $\class{K}\models\Sigma$ means
that the universal closure $\,\forall\bar{x}\,\Phi$ of each $\Phi\in\Sigma$ is
true in every structure belonging to $\class{K}$.

An \emph{atomic formula} of a first-order signature
is a formal equation $\al\thickapprox\be$ between terms or an expression $r(\al_1,\dots,\al_m)$,
where $r$ is a basic relation symbol and each $\al_j$ a term.
A \emph{strict universal Horn formula}
has the form
\begin{equation}\label{strict uh formula}
\left(\textup{\large{\&}}_{\,i<n}\,\Phi_i\right)\Longrightarrow\Phi_n,
\end{equation}
where $n\in\omega$ and $\Phi_0,\dots,\Phi_n$ are atomic formulas.
We call (\ref{strict uh formula}) a \emph{quasi-equation} if $\Phi_0,\dots,\Phi_n$
are equations.
A \emph{variety} [resp.\ \emph{quasivariety}] is the model class of a
set of atomic [resp.\ strict universal Horn] formulas.
It is said to be \emph{algebraic} if it consists of algebras $\sbA=\langle A;F\rangle$,
i.e., of structures with no indicated relation.

The class operator symbols $\III$, $\HHH$, $\SSS$, $\EEE$, $\PPP$, $\PPS$, $\PPU$ and $\RRU$ stand for
closure under isomorphic and homomorphic (surjective) images, substructures, extensions (i.e., superstructures),
direct and subdirect
products, ultraproducts and ultraroots, respectively.  Homomorphisms between similar structures are
assumed to preserve basic relations (as well as operations), but they need not reflect
the relations.  An \emph{isomorphism} is
a bijective homomorphism
whose inverse function is also a homomorphism.
An \emph{embedding} of a structure
$\sbA$ into a structure $\sbB$ is
an isomorphism from $\sbA$ onto a substructure of $\sbB$.
In a given
signature, the direct
product of the empty family is interpreted as a trivial structure.
For each class operator $\OOO$,
we abbreviate $\OOO(\{\sbA_1,\dots,\sbA_n\})$ as $\OOO(\sbA_1,\dots,\sbA_n)$.

Let $\class{K}$ be a class of similar structures.  We say that $\class{K}$ is \emph{elementary} [resp.\ \emph{universal}]
if it is the model class of a set of first-order [resp.\ universal first-order]
sentences.  This amounts to the demand that $\class{K}$ be closed under $\RRU$ [resp.\ $\SSS$],
$\III$ and $\PPU$.  Two similar structures are \emph{elementarily equivalent} if they satisfy the same first-order
sentences.
The smallest variety [resp.\ quasivariety]
containing
$\class{K}$
is $\VVV(\class{K})\seteq\HHH\SSS\PPP(\class{K})$ [resp.\
$\QQQ(\class{K})\seteq\III\SSS\PPP\PPU(\class{K})$].
The origins of these claims are discussed in \cite[Ch.\,2]{Gor98}, where proofs can also be found.
Recall that $\PPU(\class{K})\subseteq\III(\class{K})$ if $\class{K}$ is a finite set of
finite structures.

Let $\class{K}$ be a quasivariety, with $\sbA\in\class{K}$.  We say that $\sbA$ is
[\emph{finitely}] $\class{K}$-\emph{subdirectly
irreducible} if the following is true for every [finite non-empty] set $I$ and every family $\{\sbA_i:i\in I\}$ of members of
$\class{K}$\,: whenever an embedding $h\colon\sbA\mrig\prod_{I}\sbA_i$ is
subdirect
(i.e.,
$\pi_jh[A]=A_j$ for each of the projections $\pi_j\colon\prod_I\sbA_i\mrig\sbA_j$), then $\pi_i\circ h\colon\sbA\cong\sbA_i$
for some $i\in I$.  We say that $\sbA$ is $\class{K}$-\emph{simple} if it is nontrivial and every homomorphism
from $\sbA$ onto a nontrivial member of $\class{K}$ is an isomorphism.
The prefix `$\class{K}$-' is often replaced by the word `relatively' when $\class{K}$ is understood; it is
redundant when $\class{K}$ is a variety.
We denote by
$\mathsf{K}_\textit{RSI\/}$ [resp.\ $\class{K}_\textit{RFSI\,}$; $\class{K}_\textit{RS\,}$]
the class of $\mathsf{K}$-subdirectly irreducible [resp.\ finitely $\mathsf{K}$-subdirectly irreducible; $\class{K}$-simple]
members of $\mathsf{K}$.
Thus,
$\class{K}_\textit{RS}\subseteq\mathsf{K}_\textit{RSI}\subseteq\class{K}_\textit{RFSI\,}$, and $\class{K}_\textit{RSI\,}$
includes no trivial structure.  The analogue of
Birkhoff's subdirect decomposition theorem holds, i.e.,
${\mathsf{K}=\III\PPS
(\mathsf{K}_\textit{RSI\,})}$
\cite[Thm.\,3.1.1]{Gor98}.
If every $\mathsf{K}$-subdirectly irreducible member
of $\mathsf{K}$ is $\mathsf{K}$-simple, then $\mathsf{K}$ is said to be {\em relatively semisimple}.

Because the image of a homomorphism is always a substructure of the co-domain, these
definitions have the following consequence.

\begin{Fact}\label{obvious fact}
If\/ $h\colon\sbA\mrig\sbB$ is a homomorphism between members of a quasivariety,
where\/ $\sbA$ is relatively simple and\/ $\sbB$ has no trivial substructure,
then\/ $h$ is an embedding.
\end{Fact}

\section{Existential Positive Sentences}\label{eps section}

Recall that,
up to logical equivalence,
an {\em existential positive sentence\/}
is a first-order sentence of
the form
$\exists x_1 \,\dots\, \exists x_n \,\Phi$,
where $\Phi$ is a (quantifier-free) disjunction of conjunctions of atomic formulas.
Such sentences have a central place in the model theory of `positive logic' (see
\cite{PY18}, for instance).  For present purposes, their main significance derives from
Theorems~\ref{Lem:Botur} and \ref{Thm:PSC-implies-JEP} below.
They may be variable-free (and hence quantifier-free).



Given a structure $\sbA=\langle A;F;R\rangle$, with $S\subseteq A$,
let $\sbA_S=
\langle A;F\cup S_0;R\rangle$, where $S_0$ consists of the elements
of $S$, treated
as new nullary operations on $A$.
Let $\textup{Th}(\sbA)$ [resp.\ $\textup{Diag}(\sbA)$] denote the set of all [resp.\ all atomic]
first-order sentences that are true in $\sbA_A$.
%
A substructure $\sbB$ of $\sbA$ is called an \emph{elementary substructure} (and $\sbA$ an \emph{elementary extension} of $\sbB$)
if ${\sbA_B\models\textup{Th}(\sbB)}$.  In this case $\sbA$ and $\sbB$ are elementarily equivalent.
An embedding is \emph{elementary} if its image is an elementary substructure of its co-domain.
Every structure is elementarily embeddable into each of its ultrapowers.


\begin{Theorem}\label{Lem:Botur}
Let\/ $\sbA$ and\/ $\sbB$ be similar structures.  Then\/ $\sbB$ satisfies every existential positive sentence
that is true in\/ $\sbA$ iff there is a homomorphism from\/ $\sbA$ into an ultrapower of\/ $\sbB$\textup{.}
%
\end{Theorem}

\begin{proof}
($\Rig$)\,
Let $\Sigma$ be a finite subset of $\textup{Diag}(\sbA)$.  By assumption,
$\Sigma \cup \textup{Th}(\sbB)$ has a model that is an expansion of $\sbB_B$ by suitable interpretations in $B$ of the elements of $A$ occurring (as constant symbols)
in $\Sigma$.
By the Compactness Theorem, therefore,
$\textup{Diag}(\sbA) \cup \textup{Th}(\sbB)$ has a model, $\sbC$, say.  Let $\sbC^-$ be the reduct of $\sbC$ in the
signature of $\sbA,\sbB$.
Now $\sbC^-$ is isomorphic to an elementary extension of $\sbB$, because $\sbC$ is a model of $\textup{Th}(\sbB)$.
(In particular, the
negated atomic sentences in $\textup{Th}(\sbB)$ separate the elements of $B$.)
Also, there is a homomorphism from $\sbA$ into $\sbC^-$, because $\sbC$ is a model of $\textup{Diag}(\sbA)$.
As $\sbC^-$ and $\sbB$ are elementarily equivalent, some ultrapower $\sbU$ of $\sbC^-$ is isomorphic to an ultrapower of $\sbB$,
by the Keisler-Shelah Isomorphism Theorem \cite[Thm.\,6.1.15]{CK90}.  Since $\sbC^-$ is (elementarily) embeddable into $\sbU$,
there is a homomorphism from $\sbA$ into $\sbU$.

($\Leftarrow$)\, Clearly, existential positive sentences persist in homomorphic images, in
extensions and in ultraroots.
\end{proof}


\begin{Corollary}\label{Cor:model-class}
The model class of the set of existential positive sentences satisfied by a structure\/ ${\A}$ is\/ $\RRU\EEE\HHH(\sbA)$\textup{.}
\end{Corollary}

\begin{Corollary}\label{psc char}
The following demands
on a quasivariety\/ $\class{K}$ are equivalent.
\begin{enumerate}
  \item\label{psc char1}
The nontrivial members of\/ $\class{K}$ all satisfy the same existential positive sentences.

  \item\label{psc char2}
For any two nontrivial members of\/ $\class{K}$\textup{,} each can be mapped homomorphically into
an ultrapower of the other.
\end{enumerate}
\end{Corollary}

\begin{law}\label{psc def}
A quasivariety is said to be \emph{passively structurally complete} (PSC) if it satisfies the equivalent
conditions of Corollary~\textup{\ref{psc char}.}  (The terminology will be justified in Sections~\ref{sc section}
and \ref{psc section}.)
\end{law}


\section{The Joint Embedding Property}\label{jep section}


\begin{law}\label{jep defn}
A class $\mathsf{K}$ of similar structures is said to have the \textit{joint embedding property} (JEP) if, for
any two nontrivial structures ${\A},{\B}\in\mathsf{K}$, there exists ${\C}\in\mathsf{K}$ such that $\A$ and $\sbB$
can both be embedded into $\C$.
\end{law}

For quasivarieties,
the characterization of the JEP given below
was proved in \cite[Thm.\,4]{Mal66} (also see \cite[p.\,288]{Mal73} or \cite[Prop.\,2.1.19]{Gor98}).

\begin{Theorem}\label{jep char}
\textup{(Maltsev)} \,A quasivariety\/ $\class{K}$ has the JEP iff it is generated by a single structure (i.e.,
there exists\/ $\sbA\in\mathsf{K}$ such that\/ $\mathsf{K}=\QQQ(\sbA)$\textup{).}
\end{Theorem}

The following additional characterizations of the JEP for a quasivariety $\class{K}$ are
known.\footnote{\,They are specializations of conditions in \cite[Thm.\,1.2]{IG96},
where the JEP was formulated for arbitrary first-order theories, without the restriction to nontrivial
models.  It is pointed out in \cite{IG96} that their equivalence with the JEP follows from a result proved
in \cite{Zei76}.}
\begin{enumerate}
\item\label{jep2}
For each set $\mathsf{S}$ of nontrivial members of $\mathsf{K}$\textup{,} there exists a member of $\mathsf{K}$ into
which every member of $\mathsf{S}$ embeds.
\item\label{jep3}
Whenever $\Phi$ and $\Psi$ are universal sentences whose disjunction $\Phi\sqcup\Psi$ is true in all nontrivial members
of $\mathsf{K}$\textup{,} then there exists $\Xi\in\{\Phi,\Psi\}$ such that $\Xi$ is true in every nontrivial member
of $\mathsf{K}$\textup{.}
\item\label{jep4}
Whenever $\Phi$ and $\Psi$ are existential sentences, each of which is true in some nontrivial member of $\mathsf{K}$\textup{,}
then their conjunction $\Phi\aand\Psi$ is true in some nontrivial member of $\mathsf{K}$\textup{.}
\item\label{jep5}
Whenever $\Sigma$ is a set of existential sentences, each of which is true in at least one nontrivial member of $\mathsf{K}$\textup{,}
then there is a nontrivial member of $\mathsf{K}$ in which all sentences from $\Sigma$ are true.
\end{enumerate}
Easily,
(\ref{jep3}) and (\ref{jep4}) follow
from the JEP, and
(\ref{jep5}) from (\ref{jep2}).
To prove (\ref{jep2}), we apply the Compactness Theorem to
$\Sigma\cup\{\ov{\textup{Diag}}(\sbA):\sbA\in \class{S}\}$, where $\Sigma$ is a set of sentences axiomatizing $\class{K}$
and each $\ov{\textup{Diag}}(\sbA)$ is the set of atomic or negated atomic sentences that are true in $\sbA_A$.
(We arrange first that the members of $\class{S}$ are disjoint.)

Because the JEP need not persist in sub(quasi)varieties, the following result is of interest.

\begin{Theorem}\label{Thm:PSC-implies-JEP}
If a quasivariety is PSC\/ \textup{(}see Definition~\textup{\ref{psc def}),} then it has the JEP, and so do all of its subquasivarieties.
%
\end{Theorem}

\begin{proof}
Let $\sbA,\sbB$ be nontrivial members of a PSC quasivariety $\class{K}$.
Then there are homomorphisms $f \colon {\A} \mrig {\B}_{u}$ and $g \colon {\B} \mrig {\A}_{u}$, for suitable ultrapowers ${\A}_{u}$ and ${\B}_{u}$ of ${\A}$ and ${\B}$, respectively. Recall that there are (elementary) embeddings $e_{A} \colon {\A} \mrig {\A}_{u}$ and $e_{B} \colon {\B} \mrig {\B}_{u}$.
Consider the maps
\[
\langle e_{A}, f \rangle \colon {\A} \mrig {\A}_{u} \times {\B}_{u} \text{ and }\langle g, e_{B} \rangle \colon {\B} \mrig {\A}_{u} \times {\B}_{u}
\]
defined by the following rules:
for every $a \in A$ and $b \in B$,
\[
\langle e_{A}, f \rangle(a) = \langle e_{A}(a), f(a)\rangle \text{ \,and\, }\langle g, e_{B} \rangle(b) = \langle g(b), e_{B}(b)\rangle.
\]
Clearly, $\langle e_{A}, f \rangle$ and $\langle g, e_{B} \rangle$ are embeddings, so $\sbA,\sbB\in\III\SSS({\A}_{u} \times {\B}_{u})$,
and ${\A}_{u} \times {\B}_{u}  \in \QQQ({\A}, {\B})\subseteq\class{K}$.  Thus, $\class{K}$ has the JEP, as do its
sub\-quasi\-varieties, in view of the argument---or by heredity of the PSC condition.
\end{proof}

The next result allows us to restrict attention to relatively subdirectly irreducible structures when
testing a quasivariety for the JEP.

\begin{Proposition}\label{jep rsi}
Let\/ $\class{K}$ be a quasivariety, and suppose that, whenever\/ ${\sbA,\sbB\in\class{K}_\textit{RSI\,}}$\textup{,}
then there exists\/ $\sbC\in\class{K}$ such that\/ $\sbA$ and\/ $\sbB$ can both be embedded into\/ $\sbC$\textup{.}
Then\/ $\class{K}$ has the JEP.
\end{Proposition}
\begin{proof}
Let $\sbA,\sbB\in\class{K}$ be nontrivial.  Then
\[
\textup{$\sbA\in\III\PPS\{\sbA_i:i\in I\}$ \,and\,
$\sbB\in\III\PPS\{\sbB_j:j\in J\}$}
\]
for suitable $\sbA_i,\sbB_j\in\class{K}_\textit{RSI\,}$, where $I$ and $J$
are non-empty sets.  We may assume that $I\subseteq J$.  Fixing $\ell\in I$ and defining $\sbA_j=\sbA_\ell$ for all
$j\in J\backslash I$, we find that $\sbA\in\III\PPS\{\sbA_j:j\in J\}$.  By assumption, for each $j\in J$, there exists
$\sbC_j\in\class{K}$ such that $\sbA_j,\sbB_j\in\III\SSS(\sbC_j)$.  Then $\prod_J\sbA_j$ and $\prod_J\sbB_j$ both embed into
$\sbC\seteq\prod_J\sbC_j\in\class{K}$, so $\sbA,\sbB\in\III\SSS(\sbC)$.
\end{proof}

A quasivariety is said to be \emph{finitely generated} if it has the form $\QQQ(\class{K})$ for some
finite set $\class{K}$ of finite structures.  In this case, as $\QQQ=\III\PPS\SSS\PPU$,
we have $\QQQ(\class{K})_\textup{RSI}\subseteq\III\SSS(\class{K})$,
where $\SSS(\class{K})$ is again a finite set of finite structures.  Therefore,
it follows easily from Proposition~\ref{jep rsi} that
the JEP is a decidable property for finitely generated quasivarieties of finite type.

\section{Algebraic Quasivarieties and the JEP}



From now on, we confine our attention to algebraic quasivarieties.



Given an algebraic quasivariety $\class{K}$ and an algebra $\sbA$ of the same type, the $\mathsf{K}$-{\em congruences\/}
(a.k.a.\ {\em relative congruences\/})
of $\sbA$ are the congruences $\theta$ such that $\sbA/\theta\in\mathsf{K}$.  They form an
algebraic lattice ${\boldsymbol{\mathit{Con}}}_{\mathsf{K}}\sbA$, ordered by inclusion, in which
meets are intersections and the compact elements are just the finitely generated $\mathsf{K}$-congruences.
We denote by $\tcg^{\sbA}_\mathsf{K} Y$ the least $\mathsf{K}$-congruence of $\sbA$ containing a subset $Y$ of $A^2$.
When $Y=\{\langle a,b\rangle\}$, we write $\tcg^\sbA_{\mathsf{K}}(a,b)$
for the {\em principal\/} $\mathsf{K}$-congruence $\tcg^\sbA_{\mathsf{K}}Y$.
An algebra $\sbA\in\mathsf{K}$ belongs to the class $\class{K}_\textit{RSI\,}$ [resp.\ $\class{K}_\textit{RFSI\,}$;
$\class{K}_\textit{RS\,}$] iff, in ${\boldsymbol{\mathit{Con}}}_{\mathsf{K}}\sbA$, the identity relation
${\textup{id}_A\seteq\{\langle a,a\rangle:a\in A\}}$ is completely meet-irreducible [resp.\ meet-irreducible; a co-atom].
When $\mathsf{K}$ is a variety and $\sbA\in\mathsf{K}$, the congruences and $\mathsf{K}$-congruences of
$\sbA$ coincide, so the prefixes and subscripts can be dropped.

An algebra is said to be $n$-\emph{generated} (where $n$ is a cardinal) if it has a generating subset with at
most $n$ elements.  An algebraic quasivariety is said to be \emph{nontrivial} if it has a nontrivial member.
In that case, it has a relatively simple member
\cite[Thm.\,3.1.8]{Gor98}; for algebraic varieties, this was proved earlier by Magari \cite{Mag69}.
On the other hand, even when an algebra is finitely generated (i.e., $n$-generated for some $n\in\omega$), it
need not have a simple homomorphic image \cite[p.\,154]{Jon72}.
Conditions that guarantee relatively simple homomorphic images are given in the next lemma, which adapts
\cite[pp.\,153--4]{Jon72} to quasivarieties.

\begin{Lemma}\label{jonsson's other theorem}
Let\/ $\sbA$ be a nontrivial member of an algebraic quasivariety\/ $\class{K}$\textup{.}
\begin{enumerate}
\item\label{jonsson's other theorem2}
If the total relation\/ $A^2$
is compact in\/ ${\boldsymbol{\mathit{Con}}}_\mathsf{K}\sbA$\textup{,}
then\/ $\sbA$ has a relatively simple homomorphic image in\/ $\class{K}$\textup{.}

\item\label{jonsson's other theorem1}
If\/ $\sbA$ is finitely generated and of finite type, then\/ $A^2$ is compact in\/
${\boldsymbol{\mathit{Con}}}_\mathsf{K}\sbA$\textup{,} so\/ $\sbA$
has a
relatively simple homomorphic image in\/ $\class{K}$\textup{.}
\end{enumerate}
\end{Lemma}
\begin{proof}
(\ref{jonsson's other theorem2})\,  If $\bot$ is the least element of an algebraic lattice $\sbL$
and $y\in L\backslash\{\bot\}$ and $y$ is compact in $\sbL$, then $\{x\in L:y\nleqslant x\}$
has a maximal element, by Zorn's Lemma.  Setting $\sbL={\boldsymbol{\mathit{Con}}}_\mathsf{K}\sbA$ and $y=A^2$,
we conclude that, under the given assumptions, $\sbA$ has a maximal proper $\class{K}$-congruence $\theta$, whence $\sbA/\theta\in\class{K}$ is $\class{K}$-simple.

(\ref{jonsson's other theorem1})\,  Suppose $\sbA$ is generated by a finite subset $X$ of $A$.
Let $Y$ be the union of $X$ and the set of all $f(a_1,\dots,a_n)$ such that $n\in\omega$,
$f$ is a basic $n$-ary operation of $\sbA$ and $a_1,\dots,a_n\in X$.  Then
$A^2=\tcg^\sbA(Y^2)\subseteq\tcg^\sbA_\class{K}(Y^2)$.
If $\sbA$ has finite type, then $Y^2$ is finite, so $A^2=\tcg^\sbA_\class{K}(Y^2)$ is compact in
${\boldsymbol{\mathit{Con}}}_\mathsf{K}\sbA$, and the last assertion follows from (\ref{jonsson's other theorem2}).
\end{proof}

%



\begin{law}\label{kollar def}
An algebraic
quasivariety
will be called a \emph{Koll\'{a}r quasi\-variety} if each of its
nontrivial members has no trivial
subalgebra.
\end{law}

Clearly, an algebraic quasivariety $\class{K}$ is a Koll\'{a}r quasivariety if its signature includes two constant symbols
that take distinct values in every nontrivial member of $\class{K}$.  This situation is common in algebraic
logic, e.g., every quasi\-variety of
Heyting algebras
is a Koll\'{a}r quasivariety.
The result below was proved first for varieties by
Koll\'{a}r \cite{Kol79}, hence our nomenclature.
Further characterizations of Koll\'{a}r quasivarieties have been given by
Campercholi and Vaggione \cite[Prop.\,2.3]{CV12}.

\begin{Theorem}\label{kollar}
\textup{(Gorbunov \cite{Gor86}, \cite[Thm.\,2.3.16]{Gor98})}\,
An algebraic quasivariety\/ $\mathsf{K}$ is a Koll\'{a}r quasivariety iff\/ $A^2$
is compact in\/ ${\boldsymbol{\mathit{Con}}}_{\mathsf{K}}\sbA$
for every\/ ${\A} \in \mathsf{K}$\textup{.}
%
%
\end{Theorem}

\begin{Corollary}\label{kollar cor}
Every nontrivial member of a Koll\'{a}r quasivariety has a relatively simple homomorphic image.
\end{Corollary}
\begin{proof}
This follows from Theorem~\ref{kollar} and Lemma~\ref{jonsson's other theorem}(\ref{jonsson's other theorem2}).
\end{proof}

Note that an algebra
is $0$-generated iff it has a distinguished element and no proper subalgebra.
An algebra with a distinguished element has a unique $0$-generated subalgebra, which is its smallest subalgebra.
These assertions and the first two items of the next result are not always true for
structures.\footnote{\,Theorem~\ref{kollar} is generalized in \cite{Gor98} to structures having only
finitely many basic relations, using an analogue of (relative) congruences that will not be
discussed here.}

\begin{Proposition}\label{0 gen}
Let\/ $\class{K}$ be an algebraic quasivariety with the JEP.
\begin{enumerate}
  \item\label{0 gen1}
  \textup{(\cite{Los62})}
  Any two nontrivial\/ $0$-generated members of\/ $\class{K}$ are isomorphic.


  \item\label{0 gen3}
  If\/ $\class{K}$ has finite type or is a Koll\'{a}r quasivariety,
  then every nontrivial\/ $0$-generated member of\/ $\class{K}$
  is relatively simple.

  \item\label{0 gen2}
  If\/ $\class{K}$ has a constant symbol, then\/ $\class{K}$ is a Koll\'{a}r quasivariety
  or every member of\/ $\class{K}$ has a trivial subalgebra.
\end{enumerate}
\end{Proposition}
\begin{proof}
(\ref{0 gen1}) Let $\sbA,\sbB\in\class{K}$ be nontrivial and $0$-generated.  By the JEP, there exist $\sbC\in\class{K}$ and
embeddings $g\colon\sbA\mrig\sbC$
and $h\colon\sbB\mrig\sbC$.
As $g[\sbA]$ and $h[\sbB]$ are
$0$-generated substructures of $\sbC$, they
coincide,
so $h^{-1}|_{h[B]}\circ g\colon\sbA\cong\sbB$.


(\ref{0 gen3}) Let\/ $\sbA\in\class{K}$ be nontrivial and\/ $0$-generated.  The hypotheses of (\ref{0 gen3})
imply that $\sbA$ has a
homomorphic image $\sbB\in\class{K}_\textit{RS\,}$, by Corollary~\ref{kollar cor} and
Lemma~\ref{jonsson's other theorem}(\ref{jonsson's other theorem1}).  Since $\sbB$ is also $0$-generated
and nontrivial, it
is isomorphic to $\sbA$, by (\ref{0 gen1}), so $\sbA$ is relatively simple.

(\ref{0 gen2}) Let $c$ be a constant symbol of $\class{K}$.  If $\sbA\in\class{K}$ has no trivial subalgebra,
then for some basic operation symbol $f$ of $\sbA$, the atomic sentence $f(c,c,\dots,c)\thickapprox c$
(briefly, $\Phi$) is false in $\sbA$.
In that case, if
$\sbB\in\class{K}$ has a proper trivial subalgebra, then $\Phi$
is true in $\sbB$, so $\sbA$ and $\sbB$ have no common extension, contradicting the JEP.
\end{proof}

\enlargethispage{2pt}

\noindent
The assumption that $\class{K}$ has a constant symbol cannot be dropped from (\ref{0 gen2}), even when
$\class{K}$ is an algebraic variety (see Example~\ref{relevant algebras}).

\begin{Theorem}\label{Thm:simple}
Let\/ $\class{K}$ be a nontrivial Koll\'{a}r quasivariety with the JEP.
Then there is a relatively simple algebra\/ ${\A}\in\mathsf{K}$ such that\/ $\III\SSS\PPU({\A})$ includes every
relatively simple member
of\/ $\class{K}$.

Consequently, $\QQQ(\class{K}_\textit{RS\,})=\QQQ(\sbA)$\textup{,} so\/
$\QQQ(\class{K}_\textit{RS\,})$ also has the JEP.
\end{Theorem}

\begin{proof}
For any algebra $\sbB$, let $\textup{EPS}({\B})$ denote the set of existential positive sentences that are true in ${\B}$.
As $\mathsf{K}$ has the JEP, Theorem~\ref{jep char}
shows that $\class{K}=\QQQ({\C})$
for some ${\C} \in \mathsf{K}$. Since $\mathsf{K}$ is nontrivial, so is ${\C}$.
By Corollary~\ref{kollar cor},
${\C}$ has a homomorphic image ${\A}\in\class{K}_\textit{RS\,}$. Observe that
\begin{equation}\label{eps1}
\sbA\models\textup{EPS}({\C}),
\end{equation}
as $\sbA\in\HHH(\sbC)$.
We claim, moreover, that
\begin{equation}\label{eps2}
\sbC\models\textup{EPS}({\B}), \text{ for every }{\B} \in \mathsf{K}_\textit{RS\,}.
\end{equation}
Indeed, because $\class{K}=\QQQ({\C}) = \III\PPS\SSS\PPU(\sbC)$, we have
$\mathsf{K}_\textit{RS}\subseteq \mathsf{K}_\textit{RSI}\subseteq \III\SSS\PPU({\C})$, so
$\sbC\in\RRU\EEE\HHH(\sbB)$ for all $\sbB\in\class{K}_\textit{RS\,}$.
Thus, (\ref{eps2}) follows from Corollary \ref{Cor:model-class}.


Now let $\sbB\in\mathsf{K}_\textit{RS\,}$.  Then
$\sbA\models\textup{EPS}({\B})$, by (\ref{eps1}) and (\ref{eps2}),
so there is a homomorphism $h \colon {\B} \mrig \sbU$
for some ultrapower $\sbU$
of ${\A}$, by Theorem~\ref{Lem:Botur}.  Since $\sbA$ is nontrivial, so is $\sbU$.  Then $h$ is an embedding,
by Fact~\ref{obvious fact}, as $\class{K}$ is a Koll\'{a}r
quasivariety.
Thus, ${\B} \in \III\SSS\PPU({\A})$,
as claimed.

This shows that $\QQQ(\class{K}_\textit{RS\,})=\QQQ(\sbA)$, which has the JEP, by Theorem~\ref{jep char}.
\end{proof}

\begin{Corollary}\label{semisimple cor}
Let\/ $\mathsf{K}$ be a nontrivial relatively semisimple Koll\'{a}r quasi\-variety with the JEP.
Then\/ $\mathsf{K}= \QQQ({\A})$ for some relatively simple\/
${\A}\in\class{K}$\textup{.}
\end{Corollary}

\begin{proof}
This follows from Theorem \ref{Thm:simple}, as $\class{K}=\QQQ(\class{K}_\textit{RSI\,})$ and
$\class{K}_\textit{RSI}=\class{K}_\textit{RS\,}$.
\end{proof}

\begin{Corollary}
Let\/ $\mathsf{K}$ be a nontrivial Koll\'{a}r quasivariety with the JEP.
If the class
of all relatively simple members of\/ $\mathsf{K}$ is elementary, then it too has the JEP.
\end{Corollary}

\begin{proof}
Let ${\A},{\B}\in \mathsf{K}_\textit{RS\,}$.  Theorem~\ref{Thm:simple} shows that
there exist ${\C}\in\mathsf{K}_\textit{RS}$
and
${\sbU,\sbV
\in\PPU(\sbC)}$
and embeddings ${\A} \mrig \sbU$
and $\sbB\mrig\sbV$.
Because $\sbU$ and $\sbV$
are elementarily equivalent, they have a common ultrapower $\sbW$,
by the Keisler-Shelah Isomorphism Theorem, and of course $\sbW
\in\class{K}$.  Then
${\A}$ and ${\B}$ both embed into $\sbW$.
Moreover, as
$\sbW$ is elementarily equivalent to ${\C}$, and
as $\class{K}_\textit{RS}$
is elementary,
$\sbW$ is relatively simple.
\end{proof}

In view of Theorem~\ref{jep char} and Corollary~\ref{semisimple cor}, it is natural to ask whether a quasivariety
with the JEP must be generated by a relatively finitely subdirectly irreducible algebra.  Even for algebraic PSC
varieties, this is not the case, as the next example shows (also see Example~\ref{my conjecture refuted}).
Here, we require \emph{J\'onsson's Theorem} \cite{Jon67,Jon95}, which asserts that
${\VVV(\class{C})_\textit{FSI}\subseteq\HHH\SSS\PPU(\class{C})}$
for any subclass\/ $\class{C}$ of a congruence distributive algebraic variety.
Recall that every algebra with a lattice reduct generates a congruence distributive variety.


\begin{exa}\label{heyting ex}
Let $\class{K}=\VVV(\sbA,\sbB)$, where ${\A}$ and ${\B}$ are
the only two non-iso\-mor\-phic subdirectly irreducible five-element Heyting algebras.
Like every variety of Heyting algebras, $\class{K}$ is PSC and therefore has the JEP (see
Examples~\ref{Exa:PSC} below).
Suppose
$\mathsf{K}=\QQQ(\sbC)$, where $\sbC$ is
finitely subdirectly irreducible.
By J\'onsson's Theorem, $\sbC\in\HHH\SSS\PPU(\sbA,\sbB)=\HHH\SSS(\sbA,\sbB)$ (as $\sbA$ and $\sbB$ are finite),
whence $\left|C\right|\leq 5$.  Now $\sbA$ and $\sbB$ are subdirectly irreducible
members of $\QQQ(\sbC)=\III\PPS\SSS\PPU(\sbC)$, so $\sbA,\sbB\in\III\SSS\PPU(\sbC)=\III\SSS(\sbC)$
(as $\sbC$ is finite).
Since $\left|C\right|\leq\left|A\right|,\left|B\right|$,
this forces $\sbA\cong\sbC\cong\sbB$, a
contradiction.  Thus, no finitely subdirectly irreducible algebra generates $\class{K}$ as a quasivariety.
\qed
\end{exa}

An algebraic variety $\mathsf{K}$ is said to have \emph{equationally definable principal congruences} (EDPC) if there
is a finite set $\Sigma$ of pairs of $4$-ary terms in its signature
such that, whenever $\sbA\in \mathsf{K}$ and $a,b,c,d\in A$, then
\[
\langle c,d\rangle\in\tcg^\sbA
(a,b)
\text{ \,iff  }
\left(
\varphi^\sbA(a,b,c,d)=\psi^\sbA(a,b,c,d) \text{ for all  } \langle\varphi,\psi\rangle\in\Sigma
\right)\!\textup{.}
\]
In this case, $\class{K}$ is congruence distributive and has the congruence extension property (CEP),
and its class of simple members is closed under ultraproducts \cite{BP82}.

\begin{Theorem}\label{simple jep}
Let\/ $\class{K}$ be an algebraic variety with EDPC.  Then the variety\/ $\VVV(\sbA)$ has the
JEP, for every simple algebra\/ $\sbA\in\class{K}$\textup{.}
\end{Theorem}
\begin{proof}
As $\class{K}$ has EDPC, its class of simple members is closed both under $\PPU$
and (by the CEP)
under nontrivial subalgebras.
So, when
$\sbA\in\class{K}$ is simple, the nontrivial members of $\HHH\SSS\PPU(\sbA)$ belong to $\III\SSS\PPU(\sbA)$.
In this case, by J\'onsson's Theorem, $\VVV(\sbA)=\QQQ(\sbA)$, which has the JEP, by Theorem~\ref{jep char}.
\end{proof}

The JEP has a syntactic meaning in algebraic logic, which we recount below in the context of algebraic quasivarieties.
For a set $\Gamma$ of formal equations, we denote by $\mathit{Var}(\Gamma)$ the set of all variables $x$
such that $x$ occurs in at least one member of $\Gamma$.

\begin{law}\label{relevance principle}
An algebraic quasivariety $\mathsf{K}$ is said to respect the \textit{relevance principle}
if the following is true whenever
$\Gamma \cup \Delta\cup\{ \varphi \thickapprox \psi \}$
is a finite
set of equations, with
$\mathit{Var}(\Delta) \cap \mathit{Var}(\Gamma \cup \{ \varphi \thickapprox \psi \}) = \emptyset$, and $\Delta$ is consistent
over\/ $\class{K}$ (i.e., there exist terms $\al,\be$ such that
$\class{K}\not\models\left(\textup{\large{\&}}\,\Delta\right)\Longrightarrow \al\thickapprox \be$):
\[
\text{if\, }\class{K}\models\left(\textup{\large{\&}}\,(\Gamma\cup \Delta)\right) \Longrightarrow
\varphi\thickapprox \psi\text{, \,then\,
 }\class{K}\models \left(\textup{\large{\&}}\,\Gamma\right) \Longrightarrow
\varphi\thickapprox \psi.
\]
\end{law}

\begin{Theorem}\label{los suszko thm}
An algebraic quasivariety has the JEP iff it respects the relevance principle.
\end{Theorem}

This is an algebraic analogue of the \emph{\L os-Suszko Theorem} \cite[p.\,182]{LS58}, which concerns
sentential deductive systems, i.e., substitution-invariant consequence relations over terms in an algebraic signature.
(Variants of the \L os-Suszko Theorem for special families of deductive systems
are discussed in \cite{Avr14,GJKO07,KO10,Mak76,Mak95,Tok79}.)

The forward implication of Theorem~\ref{los suszko thm} follows
from item~(\ref{jep3})
after Theorem~\ref{jep char},
because a quasi-equation is logically equivalent to a disjunction (the disjuncts being its conclusion and the negations
of its premises).  Conversely, if $\sbA,\sbB$ are disjoint nontrivial members of an algebraic quasivariety $\class{K}$,
then the respective identity functions
on $A$ and $B$ extend to surjective homomorphisms $\pi_A\colon\sbF_\class{K}(A)\mrig\sbA$ and
${\pi_B\colon\sbF_\class{K}(B)\mrig\sbB}$ (where $\sbF_\class{K}(X)$ denotes a member of $\class{K}$
that is $\class{K}$-free over $X$).  In $\sbF\seteq\sbF_\class{K}(A\cup B)$, let $\theta$ be the
$\class{K}$-congruence generated by the union of the kernels of $\pi_A$ and $\pi_B$, and let $\sbC=\sbF/\theta$,
so $\sbC\in\class{K}$.  The map $a\mapsto a/\theta$ [resp.\ $b\mapsto b/\theta$] is a homomorphism from $\sbA$
[resp.\ $\sbB$] into $\sbC$.  Its injectivity follows from the relevance principle for $\class{K}$,
using the algebraicity
of the lattice $\boldsymbol{\mathit{Con}}_\class{K}\sbF$.

\begin{rema}\label{finite to infinite}
For a class $\class{K}$ of similar algebras, let $\UUU(\class{K})$ be the class of all algebras $\sbB$
such that every $\left|\!\mathit{Var}\right|$-generated subalgebra of $\sbB$ belongs to $\class{K}$.  In
general, $\UUU\III\SSS\PPP(\class{K})\subseteq\QQQ(\class{K})$, and the two need not be equal.  If $\class{K}$
is a quasivariety with the JEP, however, then there exists $\sbA\in\class{K}$ such that $\class{K}=\UUU\III\SSS\PPP(\sbA)$,
by item~(\ref{jep2}) after Theorem~\ref{jep char}, because the $\left|\!\mathit{Var}\right|$-generated members of $\class{K}$
form a set, up to isomorphism.  This means that, even if we allowed quasi-equations (over $\mathit{Var}$)
to have infinitely many premises, their validity in $\sbA$ would entail their validity throughout $\class{K}$.  In fact,
if we generalized Definition~\ref{relevance principle} to such formulas, then Theorem~\ref{los suszko thm} would
remain true.
\end{rema}

\section{Structural Completeness}\label{sc section}


In a given algebraic signature, a \emph{substitution} is an endomorphism of the absolutely free algebra
(a.k.a.\ the term algebra) generated by $\mathit{Var}$.
For an algebraic quasivariety $\mathsf{K}$ and a cardinal $m$, recall that the $\class{K}$-free $m$-generated algebra $\sbF_\class{K}(m)\in\mathsf{K}$
exists
iff $m>0$ or $\class{K}$ has a constant
symbol.
Of course, every algebraic variety $\class{K}$ is generated as such by
its free denumerably generated algebra, i.e., $\class{K}=\VVV(\sbF_\class{K}(\aleph_0))$,
but $\class{K}$ need not coincide with the quasivariety $\QQQ(\sbF_\class{K}(\aleph_0))$ (which
has the JEP, by Theorem~\ref{jep char}).

\begin{Theorem}\label{free}
\textup{(\cite[Prop.~2.3]{Ber91})}
The following conditions on an algebraic quasivariety\/ $\class{K}$ are equivalent.
\begin{enumerate}
  \item\label{free1}
  $\class{K}=\QQQ(\sbF_\class{K}(\aleph_0))$\textup.
  \item\label{free2}
  Whenever\/ $\mathsf{K}'$
  is a proper subquasivariety of\/ $\mathsf{K}$\textup{,}
  then\/ $\mathsf{K}'$ and\/ $\mathsf{K}$ generate distinct varieties, i.e.,
  $\mathbb{H}(\mathsf{K}')\subsetneq\mathbb{H}(\mathsf{K})$.
  \item\label{free3}
  For each quasi-equation\/ $\left(\varphi_1\thickapprox\psi_1\aand\dots\aand\varphi_n\thickapprox\psi_n\right)\Longrightarrow\varphi\thickapprox\psi$
  that is invalid in (some member of) $\class{K}$\textup{,} there exists a substitution $h$ such that\/
  $\class{K}\models h(\varphi_i)\thickapprox h(\psi_i)$ for\/ $i=1,\dots,n$\textup{,} but\/
  $\class{K}\not\models h(\varphi)\thickapprox h(\psi)$\textup{.}
\end{enumerate}
\end{Theorem}

\begin{law}\label{sc def}
An algebraic quasivariety $\class{K}$ is said to be \emph{structurally complete} (SC) if it satisfies the equivalent conditions of
Theorem~\ref{free}.  It is \emph{hereditarily structurally complete} (HSC) if, in addition, its subquasivarieties
are all SC.
\end{law}
In particular, an algebraic variety $\class{K}$ is SC iff each of its proper subquasivarieties generates a proper
subvariety of $\class{K}$; it is HSC iff its subquasivarieties are all varieties.  These notions have
logical origins: the algebraic
counterpart of an algebraizable deductive system $\,\vdash$ (in the sense of \cite{BP89})
is SC iff every proper extension of $\,\vdash$ has some new \emph{theorem}---as opposed to having nothing but
new rules of derivation; it is HSC iff every extension of $\,\vdash$ is axiomatic
(see for instance \cite{ORV08,PW08,Raf16}).  The terminology originates in \cite{Pog71}.

\section{Passive Structural Completeness}\label{psc section}

Every SC algebraic quasivariety $\mathsf{K}$ is PSC in the sense of Definition~\ref{psc def}.
Indeed, if $\sbA,\sbB\in\class{K}$ are nontrivial, then $\sbA$ has a homomorphic image
$\sbC\in\class{K}_\textit{RSI\,}$, while $\sbB$ is an extension of a $1$-generated homomorphic
image of $\sbF_\class{K}(\aleph_0)$, so it suffices to show that $\sbF_\class{K}(\aleph_0)$
satisfies the existential positive sentences that are true in $\sbC$.  This is indeed the
case, by Theorem~\ref{Lem:Botur},
as $\sbC\in\III\SSS\PPU(\sbF_\class{K}(\aleph_0))$ (because $\class{K}=\III\PPS\SSS\PPU(\sbF_\class{K}(\aleph_0))$, by
Theorem~\ref{free}(\ref{free1})).  Alternatively, condition~(\ref{free3}) of Theorem~\ref{free}
clearly entails the characterization of passive structural completeness in Theorem~\ref{wronski}
below.


The above argument and Theorem~\ref{Thm:PSC-implies-JEP} establish the implications
\[
\textup{HSC}\;\Longrightarrow\;
\textup{SC}\;\Longrightarrow\;
\textup{PSC}\;\Longrightarrow\;
\textup{JEP},
\]
none of which is reversible.  A variety of lattices that is SC but not HSC is exhibited in \cite[Ex.\,2.14.4]{Ber91}.
It is well known (and follows, for instance, from \cite{Min76}) that the variety of Heyting algebras is not SC, but
it is PSC (see Examples~\ref{Exa:PSC}).  An algebraic variety with the JEP that is not PSC will be pointed out in
Example~\ref{relevant algebras}.






\begin{law}
A set $\Gamma$ of equations in the signature of an algebraic quasivariety $\class{K}$ is said to be \emph{unifiable}
over $\class{K}$ if there is a substitution $h$ such that $\class{K}\models h(\varphi)\thickapprox h(\psi)$ for
every equation $\varphi\thickapprox\psi$ from $\Gamma$.
A quasi-equation
\[
(\varphi_{1} \thickapprox \psi_{1} \aand
\dots \aand
\varphi_{n} \thickapprox \psi_{n}) \Longrightarrow \varphi \thickapprox \psi
\]
in the same signature is said to be \emph{active} [resp.\ \textit{passive}] \emph{over} $\mathsf{K}$ if its set of premises
$\{\varphi_i\thickapprox\psi_i:i=1,\dots,n\}$
is [resp.\ is not] unifiable over $\class{K}$.
%
\end{law}

The next result
amplifies the logical meaning of passive structural completeness.  (It strengthens an earlier finding of
Bergman \cite[Thm.~2.7]{Ber91}.)

\begin{Theorem}\label{wronski}
\textup{(Wro\'{n}ski \cite[Fact\,2, p.\,68]{Wro09})}
An algebraic quasivariety\/ $\mathsf{K}$
is PSC iff
every quasi-equation that is passive over\/ $\mathsf{K}$ is valid in (all members of)\/ $\mathsf{K}$.
\end{Theorem}


Theorem~\ref{wronski} motivates the `passive' terminology used above, which is adapted from \cite{CM09}.
(A complementary demand, now called `active structural completeness' and analysed in \cite{DS16}, asks
that condition~(\ref{free3}) of Theorem~\ref{free} should hold for all active quasi-equations; also see
\cite{Ryb97}.  Computational aspects of [active] structural completeness are addressed in
\cite{Dyw78,MR13,Str}.)

Obviously, a quasi-equation is passive over an algebraic quasivariety $\class{K}$ iff it is passive over
the variety $\VVV(\class{K})$.  It may happen that $\class{K}$ is
PSC for the vacuous reason that no quasi-equation is passive over $\class{K}$
(as applies, for instance,
to every quasivariety of lattices).
The next result and its corollary decode this case in model-theoretic terms.
Their itemized conditions persist, of course, under varietal generation, unlike passive structural
completeness itself.


\begin{Theorem}\label{Lem:vacuos-PSC}
Let\/ $\mathsf{K}$ be an algebraic quasivariety. Then the following conditions are equivalent.
\begin{enumerate}
\item\label{vacuous1}
No quasi-equation is passive over\/ $\mathsf{K}$ (i.e., every finite set of equations in the signature
of\/ $\class{K}$ is unifiable over\/ $\class{K})$\textup{.}

\item\label{vacuous4}
$\class{K}$ is PSC and is either trivial or not a Koll\'{a}r quasivariety.

\item\label{vacuous2}
Every member of\/ $\mathsf{K}$ has an ultrapower with a trivial subalgebra.

\item\label{vacuous3}
$\sbF_{\mathsf{K}}(1)$ has an ultrapower with a trivial subalgebra.
\end{enumerate}
\end{Theorem}

\begin{proof}
(\ref{vacuous1})$\,\Rig\,$(\ref{vacuous4}):
Certainly, $\class{K}$ is PSC, by (\ref{vacuous1}) and Theorem~\ref{wronski}.
If $\sbF_\class{K}(1)$ is trivial, then every member of $\class{K}$ has a trivial subalgebra, so we may assume that
$\sbF_\class{K}(1)$ is nontrivial.

Let $\Sigma$ be the set of all existential positive sentences in the first-order
signature of $\class{K}$, and let $\{f_{1}, \dots, f_{n} \}$ be any finite set of basic operation symbols of $\mathsf{K}$.
By (\ref{vacuous1}), the equations $f_{i}(x, \dots, x) \thickapprox x$\, ($i=1,\dots,n$) are unifiable, i.e.,
there is a term $\varphi$ such that $\class{K}\models
f_i(\varphi,\dots,\varphi)\thickapprox\varphi$ for $i=1,\dots,n$.
Identifying variables, we see that
$\varphi$ may be chosen unary, whence
\begin{equation*}
\sbF_{\mathsf{K}}(1) \models
\exists x \left(x \thickapprox f_{1}(x, \dots, x) \thickapprox  \dots  \thickapprox f_{n}(x, \dots, x)\right)\!.
\end{equation*}
%
As $\{f_1,\dots,f_n\}$ was arbitrary, this implies that $\sbF_\class{K}(1)\models \Sigma$.
Let $\sbC\in\class{K}$ be trivial.
Of course, $\Sigma$ is the set of all existential positive sentences that hold in $\sbC$,
so by Theorem~\ref{Lem:Botur}, $\sbC$ can be mapped homomorphically into an ultrapower $\sbU$ of $\sbF_\class{K}(1)$,
i.e., $\sbU$ has a trivial subalgebra.  Now $\sbU$ is nontrivial (because $\sbF_\class{K}(1)\in\III\SSS(\sbU)$),
so $\class{K}$ is not a Koll\'{a}r quasivariety.
%

(\ref{vacuous4})$\,\Rig\,$(\ref{vacuous2}): Let $\sbA\in\class{K}$.  We may assume that $\class{K}$ is nontrivial (otherwise, (\ref{vacuous2}) is
immediate).  Then, by (\ref{vacuous4}), some nontrivial $\sbB\in\class{K}$ has a trivial subalgebra $\sbC$, and
$\class{K}$ is PSC, so $\sbA$ satisfies all existential positive sentences that are true
in $\sbB$.  Thus, there is a homomorphism $h\colon\sbB\mrig\sbU$ for some ultrapower $\sbU$ of $\sbA$,
by Theorem~\ref{Lem:Botur}, and $h[\sbC]$ is a trivial subalgebra of $\sbU$.

(\ref{vacuous2})$\,\Rig\,$(\ref{vacuous3}) is immediate, since $\sbF_\class{K}(1)\in\class{K}$.

(\ref{vacuous3})$\,\Rig\,$(\ref{vacuous1}):
Let $\sbU\in\PPU(\sbF_{\mathsf{K}}(1))$, where $\sbU$  has a trivial subalgebra.
Then, for any finite set $\Gamma$ of equations in the signature of $\class{K}$, the sentence
$\exists\,\ov{x}\left(\textup{\large{\&}}\,\Gamma\right)$ is true in $\sbU$, so it is true in
$\sbF_\class{K}(1)$.  Therefore, $\Gamma$
is unifiable over $\class{K}$.
%
%
%
\end{proof}

\begin{Corollary}\label{Cor:vacuous}
Let\/ $\mathsf{K}$ be an algebraic quasivariety, either of finite type or whose free\/ $1$-generated algebra is finite. Then no quasi-equation is passive over\/ $\mathsf{K}$ iff every member of\/ $\mathsf{K}$ has a trivial subalgebra.
\end{Corollary}

\begin{proof}
Sufficiency follows from Theorem~\ref{Lem:vacuos-PSC}.  Conversely, suppose that no quasi-equation is passive over $\class{K}$.
Then some ultrapower $\sbA$ of $\sbF_\class{K}(1)$ has a trivial subalgebra, again by
Theorem~\ref{Lem:vacuos-PSC}.  It clearly suffices to show that $\sbF_\class{K}(1)$ has a trivial subalgebra.
If $\sbF_\class{K}(1)$ is finite, then it is isomorphic to $\sbA$, and we are done.  If
the signature of $\class{K}$ is finite then, for its models, the property of having
a trivial subalgebra is expressed by an existential positive sentence
(which persists, of course, in ultraroots).  In that case, $\sbF_\class{K}(1)$ has a trivial subalgebra, because $\sbA$ does.
%
%
%
\end{proof}

In general, however, the ultrapowers in Theorem~\ref{Lem:vacuos-PSC}
cannot be
eliminated, because of the next example.

\begin{exa}\label{stupid algebra}
For $n\in\mathbb{N}=\{1,2,3,\dots\}$, let $f_n\colon\mathbb{N}\mrig\mathbb{N}$ be the function such that
$f_n(m)=m+n$ for $m=1,\dots,n$ and $f_n(m)=m$ whenever $n< m\in\mathbb{N}$.  Let $\sbA$ be the algebra with
universe $\mathbb{N}$, whose set of
basic operations is ${\{f_n:n\in\mathbb{N}\}}$, and let $\class{K}=\VVV(\sbA)$.  In this signature, every term that
is not a variable has the form $f_{i_1}\dots f_{i_k}(x)$ for some ${i_1,\dots,i_k\in\mathbb{N}}$.  Therefore,
since $\sbA$ generates $\class{K}$, every finite set of equations can be unified
over $\class{K}$ by substituting $f_r(x)$ for every variable, where $r$ is sufficiently large.  Thus,
no quasi-equation is passive over $\class{K}$, but $\sbA$ is a nontrivial member of $\class{K}$ that has no trivial
subalgebra.  (For each non-principal ultrafilter $\mathcal{U}$ over $\mathbb{N}$, the ultrapower
$\sbA^\mathbb{N}/\mathcal{U}$ has a trivial subuniverse, viz.\ $\{\langle 1,2,3,\dots\rangle/\mathcal{U}\}$.)
\qed
\end{exa}

Recall that an algebra $\sbA$ is said to be a {\em retract\/} of an algebra $\sbB$ if there are homomorphisms $g\colon\sbA\mrig\sbB$ and $h\colon\sbB\mrig\sbA$
such that $h\circ g$ is the identity function $\textup{id}_A$ on $A$.  This forces $g$ to be injective and $h$ surjective.

The next result identifies the PSC quasivarieties of finite type containing at least one finite nontrivial algebra.

\begin{Theorem}\label{Thm:retract}
Let\/ $\mathsf{K}$ be an algebraic quasivariety of finite type, with a finite nontrivial member.
Then the following conditions are equivalent.
\begin{enumerate}
  \item\label{retract1} $\class{K}$ is PSC.
  \item\label{retract2} The nontrivial members of\/ $\class{K}$ have a common retract.
  \item\label{retract3} Each nontrivial member of\/ $\class{K}$ can be mapped homomorphically into every member of\/ $\class{K}$\textup{.}
\end{enumerate}
In this case, the nontrivial members of\/ $\class{K}$ have a finite common retract that has no nontrivial proper
subalgebra and is either trivial or relatively simple.

Moreover, when\/ $\class{K}$ is PSC, its nontrivial members have
at most one nontrivial common retract, and they have at most one\/ $0$-generated common
retract (up to isomorphism).
%
\end{Theorem}

\begin{proof}
By assumption, $\class{K}$ has a finite nontrivial member, and
that algebra has a relatively simple (finite nontrivial) homomorphic image $\sbA\in\class{K}$.

(\ref{retract1})$\,\Rig\,$(\ref{retract2}):
Possession of a trivial subalgebra is expressible, over $\class{K}$, by an existential
positive sentence, because $\class{K}$ has finite type.  Therefore, since $\class{K}$ is PSC,
if some nontrivial member of $\class{K}$ has a trivial subalgebra, then so does every member
of $\class{K}$.  In that case, every member of $\class{K}$ has a trivial retract.


We may therefore assume that
$\class{K}$ is a Koll\'{a}r quasivariety.  In particular, $\sbA$ has no trivial subalgebra.
To complete the proof of (\ref{retract2}), we shall show that ${\A}$ is a retract of every nontrivial member of $\mathsf{K}$.

Accordingly,
let ${\B} \in \mathsf{K}$ be nontrivial, so $\sbB$ has no trivial subalgebra.  Since
${\A}$ is finite and of finite type, there is an existential positive sentence $\Phi$ such that
an algebra in the signature of $\class{K}$ satisfies $\Phi$ iff it has a subalgebra that is a
homomorphic image of $\sbA$.  As $\Phi$ is true in $\sbA$, it is true in $\sbB$, because
$\class{K}$ is PSC (and since $\sbA$ and $\sbB$ are nontrivial).  Therefore,
there is a homomorphism $g \colon {\A} \mrig {\B}$.  As ${\A}$ is relatively simple and
${\B}$ has no trivial subalgebra, $g$ is an embedding, by Fact~\ref{obvious fact}.  Moreover,
since $\mathsf{K}$ is PSC, there is a homomorphism $h$ from ${\B}$ into an ultrapower
of ${\A}$,
but ${\A}$ is finite, so $h \colon {\B} \mrig {\A}$.  Thus, $h\circ g$ is an endomorphism of $\sbA$.


Because $\sbA$ has no trivial subalgebra, the argument for the injectivity of $g$
applies equally to $h\circ g$.  Then, since $h\circ g$ is an injection from the finite
set $A$ to itself, it is surjective, i.e., $h\circ g$ is an automorphism of $\sbA$.

As the automorphism group of $\sbA$ is finite,
$(h \circ g)^{n+1}=\textup{id}_A$ for some $n\in\omega$.
Then, for the homomorphism
$k \seteq g \circ (h \circ g)^{n} \colon {\A} \mrig {\B}$,
we have $h \circ k = (h \circ g)^{n+1}=\textup{id}_A$.  Thus,
${\A}$ is a retract of ${\B}$, as claimed.

We have shown that a finite common retract $\sbA'$ of the nontrivial members of $\class{K}$ exists and can be chosen
relatively simple or trivial.  Being finite, $\sbA'$ cannot be a retract of a proper subalgebra of itself, so it has no
such nontrivial subalgebra.
In particular, if $\sbA'$ is nontrivial, then it is isomorphic to any
other nontrivial common retract of the nontrivial members of $\class{K}$.
Consequently, if $\sbA'$ is $0$-generated, then
it is isomorphic to \emph{any} other common retract of the nontrivial members of $\class{K}$,
because it is either trivial or has no trivial subalgebra.

(\ref{retract2})$\,\Rig\,$(\ref{retract3}): Let $\sbC,\sbD\in\class{K}$, where $\sbC$ is nontrivial.
We may assume that $\sbD$ is nontrivial, so there is a common retract $\sbA$ of $\sbC,\sbD$, by
(\ref{retract2}).  Then there exist a surjective homomorphism $\sbC\mrig\sbA$ and an embedding $\sbA\mrig\sbD$,
whose composition is a homomorphism $\sbC\mrig\sbD$.



(\ref{retract3})$\,\Rig\,$(\ref{retract1}):  Let $\sbC,\sbD\in\class{K}$ be nontrivial.  By (\ref{retract3}),
$\sbC$ can be mapped homomorphically into (an ultrapower of) $\sbD$,
so $\class{K}$ is PSC.
\end{proof}

\begin{Note}\label{sod finiteness}
In Theorem~\ref{Thm:retract}, the finiteness of the signature and the presence of a finite nontrivial algebra in $\class{K}$
are needed only for the implication (\ref{retract1})$\,\Rig\,$(\ref{retract2}).
\end{Note}

It follows easily from Theorem~\ref{Thm:retract}(\ref{retract3}) that passive structural completeness
is a decidable property for finitely generated algebraic quasivarieties of finite type.
Also, \ref{Thm:retract}(\ref{retract3}) amounts to the demand that each nontrivial member of $\class{K}$ is
a retract of its direct product with any member of $\class{K}$.

\begin{Corollary}\label{kollar cor 2}
Let\/ $\mathsf{K}$ be a PSC Koll\'{a}r quasivariety of finite type, with a finite nontrivial member.
Then $\class{K}$ has a unique relatively simple member (up to isomorphism), and that algebra is a
finite common retract of the nontrivial
members of\/ $\class{K}$\textup{.}
\end{Corollary}
\begin{proof}
This follows from Theorem~\ref{Thm:retract}, because a relatively simple member of a
Koll\'{a}r quasivariety is isomorphic to each of its retracts (by Fact~\ref{obvious fact}).
\end{proof}

\begin{exas}\label{Exa:PSC}
It follows from Theorem \ref{Thm:retract} that every variety consisting of groups or of Heyting algebras is PSC
(and therefore has the JEP, by Theorem~\ref{Thm:PSC-implies-JEP}).
Indeed, every nontrivial group has both a trivial retract and a subgroup with a finite nontrivial homomorphic
image, while the two-element Boolean algebra is a retract of every nontrivial Heyting algebra.
The class of all distributive lattices is a PSC variety whose nontrivial members have both
a trivial and a nontrivial common retract, the latter being the two-element lattice.
In Corollary~\ref{kollar cor 2}, we cannot drop the demand that $\class{K}$ be a Koll\'{a}r quasivariety,
as the variety of abelian groups satisfies the other hypotheses, but includes all the simple groups $\mathbb{Z}_p$
($p$ a positive prime).
\qed
\end{exas}

\begin{rema}\label{retract remark}
A $0$-generated algebra $\sbA$ is a
retract of an algebra $\sbB$ if there exist homomorphisms $g\colon\sbA\mrig\sbB$ and $h\colon\sbB\mrig\sbA$.  For in this case, every element of $A$ has the form $\al^{\sbA}(c_1,\dots,c_n)$
for some term $\al$ and some {\em distinguished\/} elements $c_i\in A$, whence $h\circ g=\textup{id}_A$, because homomorphisms preserve distinguished elements (and respect terms).
\end{rema}

\begin{notation}
\textup{For an algebraic quasivariety $\class{K}$, with $\sbA\in\class{K}$, we define
\[
\textup{Ret}(\class{K},\sbA)=
\{\sbB\in\class{K}:
\textup{$\sbB$ is trivial or $\sbA$ is a retract of $\sbB$}\}.
\]}
\end{notation}

\begin{Theorem}
\label{Thm:pre-retract}
Let\/ $\class{K}$ be an algebraic quasivariety of finite type, and\/
$\sbA\in\class{K}$
a finite\/ $0$-generated algebra.
\begin{enumerate}
  \item \label{pre-retract1}
  $\textup{Ret}(\class{K},\sbA)$ is a PSC quasivariety.


\item \label{pre-retract2.5}
  If\/ $\sbA$ is nontrivial or\/ $\class{K}$ is not a Koll\'{a}r quasivariety, then\/
  $\textup{Ret}(\class{K},\sbA)$ is a maximal PSC subquasivariety of\/ $\class{K}$\textup{.}

  \item \label{pre-retract3}
  If\/ $\class{K'}$ is a maximal PSC subquasivariety of\/ $\class{K}$\textup{,} and if\/ $\sbB'\in\class{K'}$
  is finite and nontrivial, then\/
  $\class{K'}=\textup{Ret}(\class{K},\sbA')$\textup{,} where $\sbA'$ is the\/ $0$-generated subalgebra of\/ $\sbB'$\textup{.}

  \item\label{pre-retract4}
  Every PSC subquasivariety\/
  of\/ $\class{K}$ that has a finite nontrivial member is contained in just one maximal
  PSC subquasivariety of\/ $\class{K}$\textup{.}
\end{enumerate}
\end{Theorem}
\begin{proof}
Let $\class{L}=\textup{Ret}(\class{K},\sbA)$.

(\ref{pre-retract1})\,
It suffices, by Note~\ref{sod finiteness},
to show that $\class{L}$ is a quasivariety.
As $\class{L}$
is isomorphically closed, we must show that it is closed under
$\mathbb{S}$,
$\mathbb{P}$ and $\PPU$, bearing Remark~\ref{retract remark} in mind.
If $\sbC\in\mathbb{S}(\sbB)$ and $h\colon\sbB\mrig\sbA$ is a homomorphism, then so is $h|_C\colon\sbC\mrig\sbA$,
while any embedding ${\sbA\mrig\sbB}$ maps into $\sbC$, as $\sbA$ is $0$-generated.  Thus, $\class{L}$ is closed
under $\SSS$.
Let $\{\sbB_i\colon i\in I\}$ be a subfamily of $\mathsf{L}$, where, without loss of generality, $I\neq\emptyset$.
For any $j\in I$,
the
projection $\prod_{i\in I}\sbB_i\mrig\sbB_j$ can be composed with a homomorphism
$\sbB_j\mrig\sbA$, while
$\sbA$ embeds diagonally into $\prod_{i\in I}\sbB_i$,
so $\prod_{i\in I}\sbB_i\in\mathsf{L}$.
As $\sbA$ is finite, it is isomorphic to each of its ultrapowers, so because
${\PPU\HHH(\mathsf{L'})\subseteq\HHH\PPU(\mathsf{L'})}$ for any class $\mathsf{L'}$ of similar algebras, it
follows that every ultraproduct
of $\{\sbB_i\colon i\in I\}$ can be mapped homomorphically to $\sbA$.
Also, as $\sbA$ is finite and of finite type, the attribute of having a subalgebra
isomorphic to $\sbA$ is first order-definable and therefore persists in ultraproducts.
Thus, $\mathsf{L}$ is closed under $\mathbb{P}$ and $\PPU$.



(\ref{pre-retract2.5})\, Suppose $\class{L}\subseteq\class{K'}\subseteq\class{K}$, where $\class{K'}$ is a
PSC quasivariety.  Then $\sbA\in\class{K'}$.  If $\sbA$ is nontrivial, then Theorem~\ref{Thm:retract} applies
to $\class{K'}$ (because $\sbA$ is finite) and it shows that, for every nontrivial $\sbC\in\class{K'}$, there
are homomorphisms $\sbA\mrig\sbC$ and $\sbC\mrig\sbA$ (as $\class{K'}$ is PSC).  In this case $\class{K'}\subseteq\class{L}$,
by Remark~\ref{retract remark} (as $\sbA$ is $0$-generated).  We may therefore assume that $\sbA$ is trivial.
Now suppose $\class{K}$ is not a Koll\'{a}r quasivariety.  Then $\sbA$ embeds into some nontrivial $\sbB\in\class{K}$,
whence $\sbB\in\class{L}$, and so $\sbB\in\class{K'}$.  Thus, $\class{K'}$ is not a Koll\'{a}r quasivariety.
Then $\class{K'}\subseteq\class{L}$, by Proposition~\ref{0 gen}(\ref{0 gen2}) and Theorem~\ref{Thm:PSC-implies-JEP}.




(\ref{pre-retract3})\,
Let $\class{K'},\sbB',\sbA'$ be as described.
By (\ref{pre-retract1}), it is enough to show that
${\class{K'}\subseteq\textup{Ret}(\class{K},\sbA')}$.
This will be true if every member of $\class{K'}$ has a trivial subalgebra (in which case $\sbA'$ is trivial).
We may therefore assume, by Proposition~\ref{0 gen}(\ref{0 gen2}) and Theorem~\ref{Thm:PSC-implies-JEP}, that
$\class{K'}$ is a Koll\'{a}r quasivariety (as $\class{K'}$ is PSC).
Then $\sbA'$ is nontrivial, so it is $\class{K'}$-simple, by Proposition~\ref{0 gen}(\ref{0 gen3}).  Thus,
$\class{K'}\subseteq\textup{Ret}(\class{K},\sbA')$, by Corollary~\ref{kollar cor 2}.

(\ref{pre-retract4}) follows from Theorem~\ref{Thm:retract}, together with (\ref{pre-retract1})--(\ref{pre-retract3}).
\end{proof}

A quasivariety is said to be {\em minimal\/} if it is nontrivial and has no nontrivial proper subquasivariety.
If we say that a variety is {\em minimal\/} (without
further qualification), we mean that it is nontrivial and has no nontrivial proper subvariety.  When we mean
instead that it is {\em minimal as a quasivariety}, we shall say
so explicitly, thereby avoiding ambiguity.
%
%
%
Obviously, any minimal algebraic quasivariety is HSC, and hence (P)SC.
Recall that if a (quasi)variety is finitely generated
then it is \emph{locally finite},
i.e., its finitely generated members are finite \cite[Thm.~II.10.16]{BS81}.





\begin{Theorem}\label{bergman mckenzie}
\textup{(\cite{BM90})}
\,Every locally finite congruence modular minimal algebraic variety
is also minimal as a quasivariety (and therefore HSC).
\end{Theorem}



\begin{Theorem}\label{Thm:minimal}
A relatively semisimple algebraic quasivariety\/
$\mathsf{K}$ is PSC iff it
is a minimal quasivariety or has no passive quasi-equation.
\end{Theorem}

\begin{proof}
Sufficiency is obvious.  Conversely,
let $\mathsf{K}$ be PSC and suppose that some quasi-equation is passive over $\class{K}$.
Then $\class{K}$ is a nontrivial
Koll\'{a}r quasivariety, by Theorem~\ref{Lem:vacuos-PSC}.
%
%
Let ${\A} \in \mathsf{K}$ be nontrivial.  As $\class{K}$
is relatively semisimple, its minimality will follow if we can show that
%
$\class{K}_\textit{RS}\subseteq\QQQ({\A})$,
so let ${\B} \in \mathsf{K}_\textit{RS\,}$.  Since $\sbB$ is nontrivial and $\class{K}$ is PSC,
there is a homomorphism $h$ from ${\B}$ into an ultrapower ${\C}$ of ${\A}$.  Of course, $\sbC$ is also nontrivial,
so $h$ is an embedding, by Fact~\ref{obvious fact}, because
$\class{K}$ is a Koll\'{a}r quasivariety.
Thus, $\sbB\in\III\SSS(\sbC)\subseteq\III\SSS\PPU(\sbA)\subseteq\QQQ(\sbA)$,
as required.
%
\end{proof}

The proof of Theorem~\ref{Thm:minimal} yields the following.

\begin{Corollary}
If a relatively semisimple algebraic quasivariety with a passive quasi-equation is PSC, then it is both a
Koll\'{a}r quasivariety and a minimal
quasivariety (and is therefore HSC).
\end{Corollary}

Heyting and Brouwerian algebras model intuitionistic propositional logic and its positive fragment,
respectively (see Definition~\ref{brouwerian def}).  We have noted that all varieties of Heyting algebras 
are PSC; the same applies to Brouwerian algebras, as they have trivial retracts.  Citkin has determined the 
HSC varieties of Heyting algebras \cite{Cit78} and of Brouwerian algebras \cite{Cit} (also see 
Theorem~\ref{continuum brouwerian not sc} below).  An analogous result for modal K4-algebras was proved by 
Rybakov \cite{Ryb95,Ryb97}.  Certain fragments of intuitionistic logic are modeled by HSC varieties
\cite{MW88,Pru72,Pru83,Wro86} (also see \cite{CM10}); for the case of relevance logic, see \cite{SM92} and
\cite[Sec.\,6,9]{ORV08}.
The next two sections of the present paper focus on the algebras of relevance logic (in its full signature), and their 
completeness properties.

\section{De Morgan Monoids: A Case Study}\label{dmm section}

De Morgan monoids were introduced by Dunn \cite{Dun66,MDL74}.  In the terminology of \cite{BP89},
they constitute the equivalent algebraic semantics for the principal relevance logic $\mathbf{R^t}$
of \cite{AB75}, and the quasivarieties of De Morgan monoids algebraize the extensions of
$\mathbf{R^t}$ by new axioms and/or inference rules.  There is a transparent lattice
anti-isomorphism from the subquasivarieties to the extensions, with subvarieties corresponding
to purely axiomatic extensions.  Accordingly, in \cite{MRW,MRWb} we undertook an investigation of the
lattice of varieties of De Morgan monoids.

Even when we prioritize axiomatic extensions (as relevance
logicians have tended to), the completeness conditions in Theorems~\ref{jep char} and
\ref{free}(\ref{free2}),(\ref{free3}) arise naturally and call for a consideration of
subquasivarieties as well.  The remainder of this paper therefore attempts to identify the
(quasi)varieties of De Morgan monoids that have such properties.  We describe completely
the varieties that are PSC, and those that have the JEP, and we supply some new information
concerning structural completeness.  In so doing, we are led to consider two neighbouring families
of residuated structures, viz.\ Dunn monoids and Brouwerian algebras.

\begin{law}\label{dmm def}
\textup{A
\emph{De Morgan monoid} is an algebra
$\sbA=\langle A;\bcdw,\wedge,\vee,\neg,e\rangle$
comprising a distributive lattice $\langle A;\wedge,\vee\rangle$,
a commutative monoid $\langle A;\bcdw,e\rangle$ that is \emph{square-increasing}
(i.e., $\sbA$
satisfies $x\leqslant x^2\seteq x\bcdw x$),
and a function $\neg\colon A\mrig A$,
called an {\em involution},
such that $\sbA$ satisfies $\neg\neg x\thickapprox x$ and}
\begin{equation*}
x\bcdw y\leqslant z\;\Longleftrightarrow\;x\bcdw\neg z\leqslant\neg y.
\end{equation*}
Here, $\al\leqslant\be$ abbreviates $\al\thickapprox\al\wedge\be$.
We refer to $\bcdw$ as {\em fusion}, and we define $f=\neg e$.
%
%
%
We denote by $\class{DMM}$ the class of all De Morgan monoids (which is a variety, by \cite[Thm.\,2.7]{GJKO07}).
\end{law}

For $\sbA$ as in Definition~\ref{dmm def}, fusion distributes over $\vee$, while
$\neg$ is an anti-auto\-mor\-phism of $\langle A;\wedge,\vee\rangle$, so De Morgan's laws hold.
The following facts about any De Morgan monoid $\sbA$ are known (see \cite{MRW} for sourcing).
\makeatletter
\renewcommand{\labelenumi}{\text{(\theenumi)}}
\renewcommand{\theenumi}{\Roman{enumi}}
\renewcommand{\theenumii}{\Roman{enumii}}
\renewcommand{\labelenumii}{\text{(\theenumii)}}
\renewcommand{\p@enumii}{\theenumi(\theenumii)}
\makeatother
\begin{enumerate}
\item\label{dmm one}
$\sbA$ is nontrivial iff its neutral element $e$ is not its least element.

\item\label{dmm two}
$\sbA$ is simple iff $e$ has just one strict lower bound in $\sbA$.

\item\label{dmm three}
$\sbA$ is finitely subdirectly irreducible iff $e$ is join-irreducible (or equivalently, join-prime)
in $\sbA$.

\item\label{dmm four}
$\sbA$ is subdirectly irreducible iff $e$ is completely join-irreducible in $\sbA$.

\item\label{dmm five}
If $\sbA$ has a least element $\bot$, then $a\bcdw\bot=\bot$ for all $a\in A$.

\item\label{dmm 5.5}
If $\sbA$ is finitely subdirectly irreducible and $a\in A$, then $e\leqslant a$ or $a\leqslant f$.

\item\label{dmm six}
$\sbA$ satisfies $f\leqslant e$ iff it is \emph{idempotent} (i.e., $a^2=a$ for all $a\in A$).
In this case $\sbA$ is called a \emph{Sugihara monoid}.  The \emph{odd} Sugihara monoids are
the ones in which $f=e$.

\item\label{dmm seven}
$\sbA$ satisfies $x\leqslant f^2$ iff it is \emph{anti-idempotent}, in the sense that the
variety $\VVV(\sbA)$ has no nontrivial idempotent member.

\item\label{dmm eight}
$\sbA$ satisfies $x\leqslant e$ iff it is a Boolean algebra (in which $\bcdn$ duplicates $\wedge$).

\item\label{dmm nine}
In $\sbA$, we have $f^3=f^2$.
\end{enumerate}
\makeatletter
\renewcommand{\labelenumi}{\text{(\theenumi)}}
\renewcommand{\theenumi}{\roman{enumi}}
\renewcommand{\theenumii}{\roman{enumii}}
\renewcommand{\labelenumii}{\text{(\theenumii)}}
\renewcommand{\p@enumii}{\theenumi(\theenumii)}
\makeatother


The variety $\mathsf{OSM}$ of all odd Sugihara monoids coincides with $\mathbb{V}(\sbS)$ for the algebra
$\sbS$ whose universe is the set $\mathbb{Z}$ of all integers, whose lattice order is the usual total order,
whose involution $\neg$ is additive inversion, and whose fusion is defined by
\[
a\bcdw b \, = \, \left\{ \begin{array}{ll}
                           \textup{the element of $\{a,b\}$ with greater absolute value, \,if $\left|a\right| \neq \left|b\right|$;}\\
                           a\wedge b  \textup{ \,if $\left|a\right| = \left|b\right|$.}
                                               \end{array}
                   \right.
\]
%
%
%
For each
$n\in\omega$, let $\sbS_{2n+1}$ denote the subalgebra of $\sbS$ with universe
\[
\{-n,\dots,-1,0,1,\dots,n\}.
\]
Up to isomorphism,
the algebras $\sbS_{2n+1}$ $\textup{($0<n\in\omega$)}$
are just the finitely generated subdirectly irreducible odd Sugihara monoids
(cf.\ \cite[Sec.~29.4]{AB75}).  Thus, every simple odd Sugihara monoid is isomorphic to $\sbS_3$, and
the subvariety lattice of $\class{OSM}$
is the chain
\[
\mathbb{V}(\sbS_1)\subsetneq\mathbb{V}(\sbS_3)\subsetneq\mathbb{V}(\sbS_5)\subsetneq\,\dots\,\subsetneq\mathbb{V}(\sbS_{2n+1})\subsetneq\,\dots\,\subsetneq\mathbb{V}(\sbS).
\]
%
\begin{Theorem}\label{odd sm}
\textup{(\cite{OR07,GR12})}
%
\,Every quasivariety of odd Sugihara monoids is a variety, i.e., the variety\/ $\class{OSM}$
is HSC.
\end{Theorem}


%


Infinite $1$-generated De Morgan monoids exist, but $0$-generated De Morgan monoids
are finite.  Indeed, Slaney \cite{Sla85} proved that the free $0$-generated
De Morgan monoid has exactly 3088 elements.  Its congruence lattice has just 68 elements, no two of which produce
isomorphic factor algebras \cite[Cor.\,3.6]{MRWb}.  Let $\sbA_1,\dots,\sbA_{68}$ denote the factor algebras, where
$\sbA_1$ is trivial.  By the Homomorphism Theorem, these are all of the $0$-generated De Morgan monoids, up to isomorphism.
The minimal quasivarieties of De Morgan monoids are just $\VVV(\sbS_3)$ and $\QQQ(\sbA_i)$, $i=2,\dots,68$
\,\cite[Thm.\,3.4]{MRWb}.  As passive structural completeness persists in subquasivarieties, the
next result is a characterization of the PSC quasivarieties of De Morgan monoids.

\begin{Theorem}\label{psc subq dmm}
The maximal PSC subquasivarieties of\/ $\class{DMM}$ are just the distinct classes\/ $\textup{Ret}(\class{DMM},\sbA_i)$\textup{,}
$i=1,\dots,68$\textup{,} and every nontrivial PSC quasivariety of De Morgan monoids is contained in just one of these.

Moreover,
$\textup{Ret}(\class{DMM},\sbA_1)$ is the variety of odd Sugihara monoids.  For\/ $i>1$\textup{,} each relatively simple
member of\/ $\textup{Ret}(\class{DMM},\sbA_i)$ is isomorphic to\/ $\sbA_i$\textup{.}
\end{Theorem}
\begin{proof}
A De Morgan monoid has a trivial subalgebra iff it is an odd Sugihara monoid,
so $\textup{Ret}(\class{DMM},\sbA_1)=\class{OSM}$, and
$\class{DMM}$ is not a Koll\'{a}r variety.
Therefore, $\textup{Ret}(\class{DMM},\sbA_i)$ is a maximal PSC subquasivariety of $\class{DMM}$,
for ${i=1,\dots,68}$, by Theorem~\ref{Thm:pre-retract}(\ref{pre-retract1}),\,(\ref{pre-retract2.5}).
Every maximal PSC subquasivariety $\class{K'}$ of $\class{DMM}$, other than $\class{OSM}$, has a finite
nontrivial member (viz.\ the $0$-generated subalgebra of any member of $\class{K'}\backslash
\class{OSM}$),
so $\class{K'}=\textup{Ret}(\class{DMM},\sbA_i)$ for some $i\in\{2,\dots,68\}$, by
Theorem~\ref{Thm:pre-retract}(\ref{pre-retract3}), and every nontrivial PSC subquasivariety of $\class{DMM}$ is contained
in $\textup{Ret}(\class{DMM},\sbA_i)$ for exactly one ${i\in\{1,\dots,68\}}$, by Theorem~\ref{Thm:pre-retract}(\ref{pre-retract4}).
For $i>1$, the common retract $\sbA_i$ of $\textup{Ret}(\class{DMM},\sbA_i)$ is unique (up to isomorphism) and relatively simple,
by Theorem~\ref{Thm:retract}, since it is $0$-generated and nontrivial.
\end{proof}

The PSC sub\emph{varieties} of $\class{DMM}$ are more limited.
We depict below the two-element Boolean algebra $\mathbf{2}$,
the three-element Sugihara monoid $\sbS_3$, and two
four-element
De Morgan monoids, $\sbC_4$ and $\sbD_4$.
In each case, the labeled Hasse diagram determines the structure.
Recall that $\mathbb{V}(\mathbf{2})$ is the class of all Boolean algebras.


{\tiny

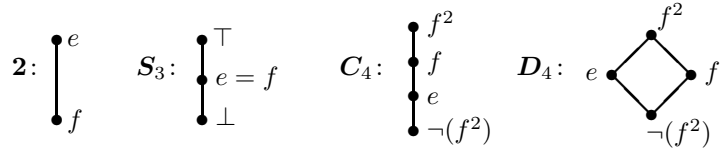
\thicklines
\begin{center}
\begin{picture}(80,60)(-28,51)

\put(-105,63){\line(0,1){30}}
\put(-105,63){\circle*{4}}
\put(-105,93){\circle*{4}}

\put(-101,91){\small $e$}
\put(-101,60){\small $f$}

\put(-122,80){\small $\mathbf{{2}\colon}$}

%
%

\put(-50,78){\circle*{4}}
\put(-50,63){\line(0,1){30}}
\put(-50,63){\circle*{4}}
\put(-50,93){\circle*{4}}

\put(-46,91){\small ${\top}$}
\put(-45,76){\small ${e}=f$}
\put(-45,61){\small ${\bot}$}

\put(-75,80){\small ${\sbS_3}\colon$}

%
%

\put(30,59){\circle*{4}}
\put(30,59){\line(0,1){39}}
\put(30,72){\circle*{4}}
\put(30,85){\circle*{4}}
\put(30,98){\circle*{4}}

\put(35,96){\small ${f^2}$}
\put(35,82){\small $f$}
\put(35,69){\small ${e}$}
\put(35,56){\small $\neg(f^2)$}

\put(2,80){\small ${\sbC_4}\colon$}

%
%

\put(120,65){\circle*{4}}
\put(135,80){\line(-1,-1){15}}
\put(135,80){\circle*{4}}
\put(105,80){\line(1,-1){15}}
\put(105,80){\circle*{4}}
\put(105,80){\line(1,1){15}}
\put(120,95){\circle*{4}}
\put(135,80){\line(-1,1){15}}

\put(122,99){\small ${f^2}$}
\put(95,78){\small ${e}$}
\put(140,78){\small $f$}
\put(118,55){\small $\neg(f^2)$}

\put(69,80){\small ${\sbD_4}\colon$}


\end{picture}\nopagebreak
\end{center}

}

\noindent
The second item of the next result is due to Slaney \cite[Thm.\,1]{Sla89}.
The first item is implicit in \cite{Sla89}; also see
\cite[Sec.\,5]{MRW}.
\begin{Theorem}\ \label{slaney}
\begin{enumerate}
  \item \label{0 gen simples}
  A De Morgan monoid is simple and\/ $0$--generated iff it is iso\-morphic to\/ $\mathbf{2}$ or to\/
  $\sbC_4$ or to\/ $\sbD_4$\textup{.}

  \item \label{slaney onto c4}
  Let\/ $h\colon\sbA\mrig\sbB$ be a homomorphism, where\/ $\sbA$ is a finitely subdirectly irreducible De Morgan monoid,
  and\/ $\sbB$ is\/ $0$--generated and nontrivial.  Then\/ $h$ is an isomorphism or\/ $\sbB\cong\sbC_4$\textup{.}
\end{enumerate}
\end{Theorem}


Every finitely generated subdirectly irreducible Sugihara monoid that is not a Boolean algebra
can be mapped homomorphically onto $\sbS_3$.  (This is well known; see for instance \cite[Sec.\,5]{MRW}.)
Consequently,
every nontrivial variety of Sugihara monoids contains $\mathbf{2}$ or $\sbS_3$.
Moreover, a variety of De Morgan monoids consists of Sugihara monoids iff it omits both $\sbC_4$ and $\sbD_4$
\cite[Thm.\,5.21]{MRW}, whence
\[
\textup{$\mathbb{V}(\mathbf{2})$\textup{,}
$\mathbb{V}(\sbS_3)$\textup{,} $\mathbb{V}(\sbC_4)$ and $\mathbb{V}(\sbD_4)$}
\]
are precisely the minimal
varieties of De Morgan monoids \cite[Thm.\,6.1]{MRW}.  (They are distinct, by J\'onsson's Theorem.)
%
%
%
%
%
Theorem~\ref{slaney}(\ref{slaney onto c4}) suggests that $\sbC_4$ has more interesting homomorphic pre-images than
the other simple $0$-generated De Morgan monoids.  We therefore make the abbreviation
\[
\mathsf{N}\seteq\textup{Ret}(\class{DMM},\sbC_4)=
\{\sbA\in\mathsf{DMM} : \,\left|A\right|=1 \textup{ \,or\, $\sbC_4$ is a retract of $\sbA$}\}.
\]
The quasivariety
$\mathsf{N}$ is
not a variety \cite[Sec.~4]{MRWb}.
It is therefore not obvious that $\class{N}$ possesses a largest subvariety,
but in fact it does.
\begin{law}\label{m defn}
We denote by $\mathsf{M}$ the variety of De Morgan monoids satisfying $e\leqslant f$ and $x\leqslant f^2$ and
\begin{equation}\label{jamie}
f^2\bcdw\neg((f\bcdw x)\wedge(f\bcdw\neg x))\thickapprox f^2\textup{.}
\end{equation}
\end{law}

\begin{Theorem}\label{m}
\textup{(\cite[Thms.~4.13,~6.8]{MRWb})}\,
$\mathsf{M}$ is the largest subvariety of\/ $\mathsf{N}$\textup{,}
and\/ $\mathsf{M}$ has\/ $2^{\aleph_0}$ distinct subvarieties.
\end{Theorem}


We can now isolate the PSC varieties of De Morgan monoids.

\begin{Theorem}\label{psc main}
Let\/ $\mathsf{K}$ be a variety of De Morgan monoids.  Then\/ $\mathsf{K}$ is PSC
iff one of the following four (mutually exclusive) conditions
holds:
\begin{enumerate}
\item\label{psc main 1}
$\mathsf{K}$ is the variety\/ $\mathbb{V}(\mathbf{2})$ of all Boolean algebras;

\item\label{psc main 2}
$\mathsf{K}=\mathbb{V}(\sbD_4)$\textup{;}

\item\label{psc main 3}
$\mathsf{K}$ consists of odd Sugihara monoids;

\item\label{psc main 4}
$\mathsf{K}$ is a nontrivial subvariety of\/ $\mathsf{M}$\textup{.}
\end{enumerate}
\end{Theorem}
\begin{proof}
By Theorem~\ref{psc subq dmm}, a nontrivial variety of De Morgan monoids is PSC iff it lies within
$\textup{Ret}(\class{DMM},\sbA_i)$ for some $i\in\{1,\dots,68\}$ (in which case $i$ is unique).
This includes
all the varieties mentioned in the present theorem, because $\mathbf{2}$, $\sbC_4$,
$\sbD_4$ and the trivial De Morgan monoid are $0$-generated and finite.  Conversely, consider a
nontrivial PSC variety ${\class{K}\subseteq \textup{Ret}(\class{DMM},\sbA_i)}$.
As $\textup{Ret}(\class{DMM},\sbA_1)=
\class{OSM}$, we may assume that $i>1$.
Theorem~\ref{psc subq dmm} also asserts that $\sbA_i$ is relatively simple in the
quasivariety $\textup{Ret}(\class{DMM},\sbA_i)$, so it is a simple member of $\class{K}$.
Therefore, $\sbA_i\in \III(\mathbf{2},\sbC_4,\sbD_4)$, by
Theorem~\ref{slaney}(\ref{0 gen simples}).  If $\sbA_i\cong\sbC_4$ then
$\class{K}\subseteq\class{M}$, by Theorem~\ref{m}, so suppose $\sbA_i\cong
\mathbf{2}$ [resp.\ $\sbA_i\cong\sbD_4$].
Let $\sbB\in\class{K}$ be subdirectly irreducible.  As $\sbA_i\in\HHH(\sbB)$,
Theorem~\ref{slaney}(\ref{slaney onto c4}) shows that $\sbB\cong\sbA_i$.  Consequently,
$\class{K}$ is $\VVV(\mathbf{2})$ [resp.\ $\VVV(\sbD_4)$].
\end{proof}

Because the four minimal varieties of De Morgan monoids are locally finite and congruence
distributive, they are minimal as quasivarieties, by Theorem~\ref{bergman mckenzie}.  In
particular, $\VVV(\mathbf{2})$ and $\VVV(\sbD_4)$ are HSC, and so is $\class{OSM}$ (as we saw in
Theorem~\ref{odd sm}).

By Theorem~\ref{psc main}, every remaining SC variety of De Morgan monoids
must be a subvariety of $\class{M}$.  By Theorem~\ref{m continuum not sc} below, $\class{M}$
and the quasivariety $\class{N}$ are not HSC.  (We conjecture that
$\class{M}$ and $\class{N}$ are not SC.)  In the lattice of subvarieties of $\class{M}$, the
unique atom $\VVV(\sbC_4)$ has just six covers, identified in \cite[Thm.\,8.10]{MRWb}.
The varietal join of those six covers is HSC \cite[Thm.\,8.13]{MRWb}.

The next result characterizes the JEP for subvarieties of $\class{DMM}$.

\begin{Theorem}\label{my conjecture}
Let\/ $\class{K}$ be a variety of De Morgan monoids.  Then\/ $\class{K}$ has the JEP iff one of the
following (mutually exclusive) conditions is met.
\begin{enumerate}
  \item\label{my conjecture1}
  $\class{K}$ is PSC.

  \item\label{my conjecture2}
  $\class{K}=\VVV(\sbA)$ for some simple De Morgan monoid\/ $\sbA$ such that\/ $\sbD_4$ is
  a proper subalgebra of\/ $\sbA$.

  \item\label{my conjecture3}
  There exist $\sbA,\sbB$ such that\/ $\class{K}=\QQQ(\sbB)$\textup{,} $\sbA$ is a simple subalgebra
  of $\sbB$\textup{,} and\/ $\sbC_4$ is a proper subalgebra of\/ $\sbA$\textup{.}
\end{enumerate}
In\/ \textup{(\ref{my conjecture3}),}
\textup{`}$\class{K}=\QQQ(\sbB)$\textup{'} can be paraphrased as\/
\textup{`}$\class{K}=\VVV(\sbB)$ and every finitely generated subdirectly irreducible member of\/ $\HHH\PPU(\sbB)$
belongs to\/ $\III\SSS\PPU(\sbB)$\textup{'.}
\end{Theorem}
\begin{proof}
Sufficiency follows from Theorems~\ref{jep char}, \ref{Thm:PSC-implies-JEP} and \ref{simple jep},
since $\class{DMM}$ has EDPC (by \cite[Thm.~3.55]{GJKO07}).

Conversely, suppose that $\class{K}$ has the JEP but is not PSC.  Then $\class{K}$ is nontrivial
and, by Theorem~\ref{psc main}, $\class{K}$ does not consist solely of Boolean algebras, nor solely of
odd Sugihara monoids.  In particular, not every member of $\class{K}$ has a trivial subalgebra.  Therefore,
$\class{K}$ is a Koll\'{a}r variety, by Proposition~\ref{0 gen}(\ref{0 gen2}), so $\sbS_3\notin\class{K}$.
As every finitely generated subdirectly irreducible Sugihara monoid that is not a Boolean algebra maps
homomorphically onto $\sbS_3$, no such algebra belongs to $\class{K}$, whence every idempotent
member of $\class{K}$ is Boolean.  Consequently, if $\class{K}$ has an idempotent member, then the $0$-generated
subalgebras of its nontrivial members are all isomorphic to $\mathbf{2}$, by Proposition~\ref{0 gen}(\ref{0 gen1}).
In that case, $\class{K}$ consists of idempotent algebras, by (\ref{dmm six}), and so coincides with
$\VVV(\mathbf{2})$, a contradiction.  This shows that $\class{K}$ has no nontrivial idempotent member.

Being nontrivial, $\class{K}$
therefore includes $\sbC_4$ or $\sbD_4$, so $\III(\sbC_4)$ or $\III(\sbD_4)$ is the class of all
$0$-generated nontrivial members of $\class{K}$, by Proposition~\ref{0 gen}(\ref{0 gen1}).
Also, by (\ref{dmm seven}),
$\class{K}$ satisfies $x\leqslant f^2$ (and hence $\neg(f^2)\leqslant x$ as well).  On the other hand,
$\class{K}\not\subseteq\class{M}$ and $\class{K}\neq\VVV(\sbD_4)$, by Theorem~\ref{psc main}, as $\class{K}$
is not PSC.

By Theorem~\ref{Thm:simple}, there is a simple De Morgan monoid $\sbA\in\class{K}$ such that
\begin{equation}\label{a}
\textup{every simple member of $\class{K}$
belongs to $\III\SSS\PPU(\sbA)$.}
\end{equation}
By Theorem~\ref{jep char}, there exists $\sbE\in\class{K}$ such that $\class{K}=\QQQ(\sbE)=\III\PPS\SSS\PPU(\sbE)$,
whence $\sbA\in\III\SSS\PPU(\sbE)$ (as $\sbA$ is simple).  Choose $\sbB\in\III\PPU(\sbE)$ with
$\sbA\in\SSS(\sbB)$.  As $\sbE\in\III\SSS(\sbB)$, we have $\class{K}=\QQQ(\sbB)$.

Suppose first that $\III(\sbC_4)$ is the class of $0$-generated nontrivial members of $\class{K}$.
As $\class{K}$ is a Koll\'{a}r variety, the $0$-generated subalgebra of $\sbA$ is nontrivial, so it
can be identified with $\sbC_4$.  If $\sbA=\sbC_4$, then every simple member
of $\class{K}$ is isomorphic to $\sbC_4$, by (\ref{a}), so $\sbC_4$ is a retract of every
nontrivial member of $\class{K}$ (by Corollary~\ref{kollar cor} and Remark~\ref{retract remark}),
i.e., $\class{K}\subseteq\class{M}$, a contradiction.
This shows that $\sbC_4$ is a proper subalgebra of $\sbA$, so (\ref{my conjecture3}) holds.

We may now assume that $\III(\sbD_4)$ is the class of $0$-generated nontrivial members of $\class{K}$.
Let $\sbG$ be any subdirectly irreducible member of $\class{K}$.
Again, since $\class{K}$ is a Koll\'{a}r variety, the $0$-generated subalgebras of $\sbA,\sbG$ are nontrivial,
so we may assume that $\sbD_4\in\SSS(\sbA)\cap\SSS(\sbG)$.  Therefore, as $\sbD_4$ satisfies
$e\wedge f\thickapprox\neg(f^2)$, so does $\sbG$.  Consequently, as $\neg(f^2)$ is the least element of $\sbG$,
it follows from (\ref{dmm 5.5}) that $\neg(f^2)$ is the greatest strict lower bound of $e$ in $\sbG$, whence
$\sbG$ is simple.  This shows that $\class{K}$
is a semisimple variety, so $\class{K}=\QQQ(\sbA)$, by (\ref{a}).
Since $\class{K}\neq\VVV(\sbD_4)=\QQQ(\sbD_4)$, we must have $\sbA\neq\sbD_4$, and so
(\ref{my conjecture2}) holds.

Note that (\ref{my conjecture1}) precludes both (\ref{my conjecture2}) and (\ref{my conjecture3}),
by Theorem~\ref{Thm:retract}, because each of $\sbC_4,\sbD_4$ has no retract other than its isomorphic images,
and cannot be a retract of a strictly larger simple algebra.  Also, (\ref{my conjecture2}) precludes
(\ref{my conjecture3}), by Proposition~\ref{0 gen}(\ref{0 gen1}), as $\sbC_4$ and $\sbD_4$ are both
$0$-generated and nontrivial.

As every variety is generated as such by its finitely generated subdirectly irreducible members, the paraphrase
in the last claim is justified by J\'{o}nsson's Theorem and the CEP for De Morgan monoids
(which implies that ${\HHH\SSS(\sbP)\subseteq\SSS\HHH(\sbP)}$ for all $\sbP\in\class{DMM}$).
\end{proof}

\begin{Corollary}
A variety of Sugihara monoids has the JEP iff it is PSC.
\end{Corollary}

\begin{Corollary}
In the lattice of varieties of De Morgan monoids, all but four of the join-irreducible covers of
atoms have the JEP.
\end{Corollary}
\begin{proof}
The join-irreducible covers of the four atoms are described in \cite[Thm.\,7.2]{MRWb}.
With four exceptions, each has the JEP, as it is either a subvariety of $\class{M}$ or of $\class{OSM}$
(and is thus PSC) or has the form $\VVV(\sbA)$ for a simple algebra $\sbA$.
The exceptions are covers of $\VVV(\sbC_4)$ that lack the JEP,
by Proposition~\ref{0 gen}(\ref{0 gen1}), because they are the varietal closures of
$0$-generated algebras $\sbC_5,\dots,\sbC_8$ (respectively), each of which has more elements than $\sbC_4$.
\end{proof}

In Theorem~\ref{my conjecture}(\ref{my conjecture3}), it can happen that $\class{K}$ is not generated, even as a variety, by one finitely
subdirectly irreducible algebra.  This remains the case when $\class{K}=\QQQ(\sbB)$ for some \emph{finite} $\sbB$.
These claims will be justified in Example~\ref{my conjecture refuted}.

\begin{exa}\label{relevant algebras}
A \emph{relevant algebra} is an $e$-free subreduct of a De Morgan monoid (i.e., a subalgebra of the reduct $\langle A;\bcdw,\wedge,\vee,\neg\rangle$ of some $\sbA\in\class{DMM}$).
These algebras form a variety $\class{RA}$, algebraizing the fragment $\mathbf{R}$ of $\mathbf{R^t}$ that lacks
the so-called
Ackermann constants.  A finite equational basis
for $\mathsf{RA}$ is
given in \cite{FR90} (also see \cite{Dzi83}, \cite[Cor.\,4.11]{HR07} and \cite[Sec.\,7]{MRW}).
Boolean algebras may be regarded as relevant algebras, since they satisfy ${e\thickapprox x\vee\neg x}$.

The variety $\class{RA}$ respects the relevance principle of Definition~\ref{relevance principle}, by
\cite[Thm.\,6]{Mak76} (since every finite set of equations in its signature is consistent over $\class{RA}$,
as follows from a consideration of the locally finite $e$-free reduct of the odd Sugihara monoid $\sbS$).  Therefore,
$\class{RA}$ has the JEP, by Theorem~\ref{los suszko thm}.  In other words (by Theorem~\ref{jep char}),
$\class{RA}=\QQQ(\sbA)$ for some
$\sbA$ (cf.\ \cite[Thm.\,5]{Tok79}).  This $\sbA$ is not $\sbF_\class{RA}(\aleph_0)$, as
$\class{RA}$ is not (even passively)
structurally complete.  In fact, $\class{RA}$ has no nontrivial PSC subvariety, other than $\VVV(\mathbf{2})$
\cite[Thm.\,6]{RS16}.

This contrasts with the fact that
$\class{DMM}$ lacks the JEP (by Proposition~\ref{0 gen}(\ref{0 gen1}), as it has non-isomorphic $0$-generated nontrivial
members).
Because $\mathbf{2}$ has no trivial subalgebra,
while the $e$-free reduct of $\sbS_3$ has a trivial subalgebra (and so belongs to no Koll\'{a}r quasivariety),
$\class{RA}$ would violate Proposition~\ref{0 gen}(\ref{0 gen2}) if we dropped the demand there for a constant symbol in the
signature.
\qed
\end{exa}

We have seen that $\class{M}$ is PSC and contains all the SC subvarieties of $\class{DMM}$ not explicitly
identified in Theorem~\ref{psc main}.  We shall show, however, that $\class{M}$ has $2^{\aleph_0}$ structurally incomplete
subvarieties.  For this we need to consider algebras called Dunn monoids, and
a construction known as reflection.



\section{Dunn Monoids and Reflections}\label{dunn monoids and reflections}

With respect to the derived operation
$x\rig y\seteq\neg(x\bcdw\neg y)$,
every De Morgan monoid satisfies
$\neg x\thickapprox x\rig f$ and
\begin{equation}\label{residuation}
x\bcdw y\leqslant z\;\Longleftrightarrow\;y\leqslant x\rig z \quad\textup{(the law of residuation).}
\end{equation}

\begin{law}\label{dm def}
An algebra
$\sbA=\langle A;\bcdw,\rig,\wedge,\vee,e\rangle$
is called a \emph{Dunn monoid}
if $\langle A;\wedge,\vee\rangle$ is a distributive lattice,
$\langle A;\bcdw,e\rangle$ is a square-increasing commutative monoid
and
$\rig$ is a binary operation---called \emph{residuation}---such
that $\sbA$ satisfies
(\ref{residuation}).
\end{law}
%
%
%
Dunn monoids form a variety, again by \cite[Thm.\,2.7]{GJKO07}.  This variety is PSC, because $\{e\}$
constitutes a subalgebra of any Dunn monoid.
%
Clearly, up to term equivalence, every De Morgan monoid has a reduct that is
a Dunn monoid.
Conversely, as recounted below, each Dunn monoid can be embedded into (the $\bcdw,\rig,\wedge,\vee,e$
reduct of) a De Morgan monoid.
Properties (\ref{dmm one})--(\ref{dmm five})
of De Morgan monoids remain true for Dunn monoids.
\begin{law}\label{reflection def}
\textup{(cf.\ Meyer \cite{Mey73})}\,
Given a Dunn monoid $\sbA$
and a disjoint copy ${A'=\{a':a\in A\}}$ of $A$, let $\bot,\top$ be distinct non-elements of $A\cup A'$.
By the \emph{reflection} $\textup{R}(\sbA)$ of $\sbA$, we mean the De Morgan monoid with universe
${\textup{R}(A)=A\cup A'\cup\{\bot,\top\}}$ such that
$\sbA$ is a subalgebra of the Dunn monoid reduct of $\textup{R}(\sbA)$ and, for all $a,b\in A$ and $x\in \textup{R}(A)$,
\begin{align*}
& x\bcdw\bot=\bot < a < b' < \top = a'\bcdw b', \textup{ and if $x\neq\bot$, then }x\bcdw\top=\top;\\
& a\bcdw b' = (a\rig b)';\\
& \neg a = a' \textup{ and } \neg(a')=a \textup{ and } \neg\bot=\top \textup{ and } \neg\top = \bot.
\end{align*}
\end{law}
In Definition~\ref{reflection def}, since $f=e'$, we have $\top=f^2$ and $\bot=\neg(f^2)$, so $\bot,\top$
belong to every subalgebra of $\textup{R}(\sbA)$.  The subalgebra on $\{\bot,e,e',\top\}$ is isomorphic to $\sbC_4$.
The De Morgan monoid $\textup{R}(\sbA)$ satisfies $e\leqslant f$ and $x\leqslant f^2$ and (\ref{jamie}),
so it belongs to $\class{M}$.
The construction of $\textup{R}(\sbA)$ from $\sbA$ preserves and reflects [finite] subdirect
irreducibility---except that $\textup{R}(\sbA)$ is subdirectly irreducible when $\sbA$ is trivial.
The interaction between reflections and the operators $\SSS,\HHH,\PPU$ is illuminated in the next lemma.

%

\begin{Lemma}
\label{hs lem}
\textup{(\cite[Lem.\,6.5]{MRWb})}
Let\/ $\sbA$ be a Dunn monoid.
\begin{enumerate}
\item\label{s}
If\/ $\sbB$ is a subalgebra of $\sbA$\textup{,} then\/
$B\cup\{b':
b\in B\}\cup\{\bot,\top\}$
is the universe of a subalgebra of\/ $\textup{R}(\sbA)$ that is isomorphic to $\textup{R}(\sbB)$\textup{,} and every
subalgebra of\/ $\textup{R}(\sbA)$ arises in this way from a subalgebra of\/ $\sbA$\textup{.}

\smallskip
\item\label{h}
If\/ $\theta$ is a congruence of\/ $\sbA$\textup{,} then
\[
\textup{\quad\quad $\textup{R}(\theta)\seteq\theta\cup\{\langle a',b'\rangle:
\langle a,b\rangle\in \theta\}\cup\{\langle \bot,\bot\rangle,\,\langle \top,\top\rangle\}$}
\]
is a congruence of\/ $\textup{R}(\sbA)$\textup{,} and\/ $\textup{R}(\sbA)/\textup{R}(\theta)\cong\,\textup{R}(\sbA/\theta)$\textup{.}  Also, every proper congruence of\/
$\textup{R}(\sbA)$ has the form\/ $\textup{R}(\theta)$ for some\/ $\theta\in{{\mathit{Con}}}\,\sbA$\textup{.}

\smallskip

\item\label{pu}
If\/ $\{\sbA_i:
i\in I\}$ is a family of Dunn monoids
and\/ $\mathcal{U}$ is an ultrafilter over\/ $I$\textup{,} then\/
$\prod_{i\in I}\textup{R}(\sbA_i)/\mathcal{U}\,\cong\textup{R}\!\left(\prod_{i\in I}\sbA_i/\mathcal{U}\right)$\textup{.}
\end{enumerate}
\end{Lemma}

\begin{law}\label{reflection def2}
The \emph{reflection}
of a variety $\mathsf{K}$ of Dunn monoids
is the
subvariety
$\mathbb{R}(\mathsf{K})\seteq\mathbb{V}\{\textup{R}(\sbA):\sbA\in\mathsf{K}\}$
of $\class{M}$.
\end{law}
As a function from the lattice of varieties of Dunn monoids to the subvariety lattice
of $\mathsf{M}$, the operator $\mathbb{R}$
is obviously isotone.  It is also order-reflecting, and therefore injective \cite[Lem.\,6.7]{MRWb}.


\begin{Theorem}\label{r structural incompleteness}
Let\/ $\mathsf{K}$ be a variety of Dunn monoids.  If\/ $\mathbb{R}(\mathsf{K})$ is structurally complete, then so is\/ $\mathsf{K}$
(i.e., $\mathbb{R}$ preserves structural\/ \textup{in}completeness).
\end{Theorem}
\begin{proof}
Suppose $\mathsf{K}$ is not SC,
so $\mathsf{K}=\mathbb{H}(\mathsf{L})$ for some
quasivariety $\mathsf{L}\subsetneq
\mathsf{K}$.  Now $\mathsf{L}^\dagger\seteq\mathbb{I}\{\textup{R}(\sbB):\sbB\in\mathsf{L}\}$ is closed under $\mathbb{S}$ and $\mathbb{P}_\mathbb{U}$,
by Lemma~\ref{hs lem}(\ref{s}),(\ref{pu}),
so $\mathbb{Q}(\mathsf{L}^\dagger)=\mathbb{IP}_\mathbb{S}(\mathsf{L}^\dagger)\subseteq\mathbb{R}(\mathsf{K})$.  As $\mathsf{L}\subsetneq\mathsf{K}$,
and because all quasivarieties are closed under subdirect products, there is an algebra $\sbA\in\mathsf{K}_\textit{SI\,}\backslash
\mathsf{L}$.
Then $\textup{R}(\sbA)$ belongs to $\mathbb{R}(\mathsf{K})$ and is subdirectly irreducible.  \,So, if ${\textup{R}(\sbA)\in\mathbb{Q}(\mathsf{L}^\dagger)}$,
then
$\textup{R}(\sbA)\cong\textup{R}(\sbB)$ for some $\sbB\in\mathsf{L}$, whence $\sbA\cong\sbB$, contradicting the fact that $\sbA\notin\mathsf{L}$.
Therefore, $\textup{R}(\sbA)\notin\mathbb{Q}(\mathsf{L}^\dagger)$, and so
${\mathbb{Q}(\mathsf{L}^\dagger)\neq\mathbb{R}(\mathsf{K})}$.

We claim that $\mathbb{R}(\mathsf{K})=\mathbb{V}(\mathsf{L}^\dagger)$.  To see this, let $\sbC\in\mathbb{R}(\mathsf{K})_\textit{SI\,}$.  \,By
J\'{o}nsson's Theorem and Lemma~\ref{hs lem}, $\sbC\cong\textup{R}(\sbD)$ for some $\sbD\in\mathsf{K}$.  As ${\mathsf{K}=\mathbb{H}(\mathsf{L})}$, we may assume
that $\sbD=\sbE/\theta$ for some $\sbE\in\mathsf{L}$ and some ${\theta\in{{\mathit{Con}}}\,\sbE}$.  Then
$\sbC\cong\textup{R}(\sbE/\theta)\cong\textup{R}(\sbE)/\textup{R}(\theta)$,
by Lemma~\ref{hs lem}(\ref{h}), whence
${\sbC\in\mathbb{H}(\mathsf{L}^\dagger)\subseteq\mathbb{V}(\mathsf{L}^\dagger)}$.  This vindicates the claim.

In summary, $\mathbb{Q}(\mathsf{L}^\dagger)$ is a proper subquasivariety of $\mathbb{R}(\mathsf{K})$ that fails to generate a \emph{proper}
subvariety of $\mathbb{R}(\mathsf{K})$, so $\mathbb{R}(\mathsf{K})$ is not SC.
\end{proof}

A Dunn (or De Morgan) monoid is said to be \emph{semilinear} if it is isomorphic to a subdirect product
of totally ordered algebras.  The semilinear Dunn monoids form a structurally incomplete variety
\cite[Thm.~9.4]{ORV08}.  Its reflection is just the variety $\mathsf{SLM}$
of semilinear members of $\mathsf{M}$, by Lemma~\ref{hs lem} and \cite[Cor.~5.8]{MRWb},
so $\mathsf{SLM}$ is not structurally complete either, by Theorem~\ref{r structural incompleteness}.
This confirms that $\class{M}$ is not HSC (and likewise $\class{N})$, but we can say more:

\begin{Theorem}\label{m continuum not sc}
The variety\/ $\class{M}$ has\/ $2^{\aleph_0}$ structurally incomplete subvarieties.
\end{Theorem}
\begin{proof}
Since the operator $\mathbb{R}$ is injective, it suffices, by Theorem~\ref{r structural incompleteness},
to exhibit $2^{\aleph_0}$ structurally incomplete varieties of Dunn monoids.  The existence of such a
family will be proved in Theorem~\ref{continuum brouwerian not sc} below.
\end{proof}

We conclude this section with an example illustrating Theorem~\ref{my conjecture}(\ref{my conjecture3}).

\begin{exa}\label{my conjecture refuted}
The Dunn monoid reduct of a De Morgan monoid $\sbA$ shall be denoted by $\sbA^+$.
We then denote by $X(\sbA)$ the De Morgan monoid that extends $\textup{R}(\sbA^+)$
by just one element $x$, where $a<x<b'$ for all $a,b\in A$, and $x\bcdw\neg(f^2)=\neg(f^2)$
and $x=\neg x=x\bcdw c$ and $x\bcdw d=f^2$ whenever $\neg(f^2)<c\leqslant x<d\leqslant f^2$.
(It is easily checked that this $X(\sbA)$ is indeed a De Morgan monoid,
with ${\textup{R}(\sbA^+)\in\SSS(X(\sbA))}$.)

Let $\class{K}=\VVV(X(\mathbf{2}\times\sbS_3))$.
As $\class{K}$ is
generated by one finite algebra, its finitely subdirectly irreducible
members are finite and can be computed mechanically, by J\'{o}nsson's Theorem.  None of them
has the property that its $\HHH\SSS$-closure contains
all the others, but all of them embed into $X(\mathbf{2}\times\sbS_3)$ (excepting the trivial algebra).  Therefore,
$\class{K}$ is not generated as a variety by a single (finitely) subdirectly irreducible algebra, but $\class{K}=\QQQ(X(\mathbf{2}\times\sbS_3))$, so $\class{K}$ has the JEP, by Theorem~\ref{jep char}.
For the trivial De Morgan monoid $\sbE$, the
five-element simple algebra $X(\sbE)$ belongs to $\class{K}$ and has $\sbC_4$ as its smallest subalgebra.
As cases (\ref{my conjecture1}) and (\ref{my conjecture3}) of Theorem~\ref{my conjecture} are mutually
exclusive, $\class{K}$ is not PSC (in contrast with Example~\ref{heyting ex}).
\qed
\end{exa}

\section{Brouwerian Algebras}\label{brouwerian section}

\begin{law}\label{brouwerian def}
A Dunn monoid is called a \emph{Brouwerian algebra} if it satisfies
$x\bcdw y\thickapprox x\wedge y$ (or equivalently, $x\leqslant e$),
in which case it
is identified with its $\rig,\wedge,\vee,e$ reduct.
\end{law}


A Heyting algebra is therefore just a Brouwerian algebra with a distinguished least element.
A well-known duality exists between Heyting algebras and `Esakia spaces' \cite{Esa74,Esa85}.
It entails a duality between Brouwerian algebras and `pointed Esakia spaces' (see \cite[Sec.\,3]{BMR17}, for instance).
Here we require only a topology-free version of the latter duality, explained briefly below.

In an indicated partially ordered set (`poset') $\sbX=\langle X;\leqslant\rangle$,
we define
\[
\textup{$\up{x}=\{y\in X\colon x\leqslant y\}$
\ and \ $\up U=\mbox{\scriptsize $\bigcup$}_{u\in U}\up{u}$,}
\]
for $U\cup\{x\}\subseteq X$, and if $U=\up U$, we call $U$ an {\em up-set\/} of $\sbX$.
We define $\down x$ and $\down U$
dually.  (Where ambiguity is a danger, we write $\up{x}$ as $\up^{\sbX\!}{x}$, etc.)
We call $\sbX$ a \emph{dominated} [resp.\ \emph{bounded}] poset if it
has a greatest [resp.\ a greatest and a least] element.
When $\sbX$ is dominated, the set $\textit{Up\/}(\sbX)$ of all non-empty up-sets of $\sbX$ is closed
under intersections and if, for $U,V\in\textit{Up\/}(\sbX)$, we define
\[
U\rig V\seteq X\backslash\!\downarrow\!(U\backslash V)\;\;\left(=\mbox{$\bigcup$}\,\{W\in\textit{Up\/}(\sbX):W\cap U\subseteq V\}\right)
\]
then $\sbX^*\seteq\langle \textit{Up\/}(\sbX);\rig,\cap,\cup,X\rangle$
is a Brouwerian algebra.  Note that $\sbX^*$ is subdirectly irreducible iff $\sbX$ is bounded and not a singleton.

A function $g\colon\sbX\mrig\sbY$ between dominated posets is called a $\emph{$p$-morphism}$ if it is isotone and
\begin{equation}\label{p-morphism}
\textup{whenever $g(x)\leqslant y\in Y$, then $y=g(z)$ for some $z\in\up{x}$.}
\end{equation}
In this case, $\up{g(x)}=g[\up{x}]$ for all $x\in X$ (whence $g[X]\in\textit{Up\/}(\sbY)$), $g$ preserves top elements,
and there is a
homomorphism $g^*\colon\sbY^*\mrig\sbX^*$, defined by $V\mapsto g^{-1}[V]\seteq
\{x\in X:g(x)\in V\}$ ($V\in\textit{Up\/}(\sbY)$).
In particular, each $U\in\textit{Up\/}(\sbX)$ is the universe of
a dominated subposet $\sbU$ of $\sbX$, and since the inclusion map $i\colon\sbU\mrig\sbX$ is a $p$-morphism, the
function $i^*\colon V\mapsto U\cap V$ ($V\in\textit{Up\/}(\sbX)$) is a
homomorphism $\sbX^*\mrig\sbU^*$.
As it is obviously surjective, ${\sbU^*\in\HHH(\sbX^*)}$ for all $U\in\textit{Up\/}(\sbX)$.

For a Brouwerian algebra $\sbA$, we denote by $\Pr(\sbA)$ the set of all prime filters of the lattice
$\langle A;\wedge,\vee\rangle$, \emph{including} $A$ itself.  Thus, $\Pr(\sbA)$ consists of the non-empty up-sets
$P$ of the poset reduct $\langle A;\leqslant\rangle$ of $\sbA$ such that $P$ is closed under $\wedge$ and
$A\backslash P$ is closed under $\vee$.  The dominated poset $\langle\Pr(\sbA);\subseteq\rangle$ is
abbreviated as $\sbA_*$.  It is bounded iff $\sbA$ is finitely subdirectly irreducible.  If $h\colon\sbA\mrig\sbB$
is a homomorphism between Brouwerian algebras, then there is a $p$-morphism $h_*\colon\sbB_*\mrig\sbA_*$, defined by
$Q\mapsto h^{-1}[Q]$ ($Q\in\Pr(\sbB)$).

For each Brouwerian algebra $\sbA$, there is an embedding
$\sbA\mrig{\sbA_*}^*$, defined by $a\mapsto\{P\in\Pr(\sbA):a\in P\}$; it is an
isomorphism if $\sbA$ is finite.  For each dominated poset $\sbX$, there is an injective $p$-morphism $\sbX\mrig{\sbX^*}_{*\,}$,
defined by $x\mapsto\{U\in\textit{Up\/}(\sbX):x\in U\}$; it is bijective when $\sbX$
is finite, in which case its inverse is also a $p$-morphism.

The functor $\sbA\mapsto\sbA_{*}$\,; $h\mapsto h_*$ defines a duality from the category $\class{FBA}$ of finite Brouwerian
algebras (and their homomorphisms) to the category $\class{FDP}$ of finite dominated posets (and their $p$-morphisms),
i.e., it defines a category equivalence from $\class{FBA}$ to the opposite category of $\class{FDP}$.  A reverse
functor is given by $\sbX\mapsto\sbX^*$\,; $g\mapsto g^*$.

In particular, every finite dominated poset is isomorphic to
$\sbA_*$ for some finite Brouwerian algebra $\sbA$, and if $\sbA,\sbB$ are finite Brouwerian algebras, then the rule $h\mapsto h_*$
defines a bijection from the set of all homomorphisms $\sbA\mrig\sbB$ to the set of all $p$-morphisms $\sbB_*\mrig\sbA_*$.

One feature of the above duality is the following (cf.\ \cite[Lem.\,3.4(ii)]{BMR17}).

\begin{Lemma}\label{1-1 and onto}
A homomorphism\/ $h$ between finite Brouwerian algebras
is surjective iff\/ $h_*$ is injective.  Also, $h$ is injective iff\/ $h_*$ is surjective.

Consequently, a\/ $p$-morphism\/ $g$ between finite dominated posets is surjective\/ [resp.\ injective] iff\/
$g^*$ is injective [resp.\ surjective].
\end{Lemma}

\begin{law}\label{depth defn}
In a dominated poset $\sbX$, the \emph{depth} of
an element $x$
is
the largest non-negative integer $n$ (if it exists) such that the
subposet $\up{x}$ contains a
chain of cardinality $n+1$.
Thus, the greatest element of $\sbX$  has depth $0$.
If $n$ is minimal such that all elements of $\sbX$ have depth at most $n$, then
$\sbX$ itself is said to have \emph{depth} $n$.
A variety $\class{K}$ of Brouwerian algebras is said to have \emph{depth} $n$ if
$\sbA_*$
has depth at most $n$ for every $\A \in \class{K}$.
\end{law}
By (\ref{p-morphism}), for any
$p$-morphism
$g\colon\sbX\mrig\sbY$ between dominated posets,
\begin{equation}\label{depth preservation}
\textup{if $x\in X$ has depth $n\in\omega$, then $g(x)$ has depth at most $n$.}
\end{equation}

Each dominated poset $\sbX$ is an up-set of a dominated poset $\widehat{\boldsymbol{X}}$,
which differs from $\sbX$ only as follows: whenever $a,b$ are distinct elements of depth $2$
in $\sbX$, then $\widehat{\sbX}$ has a (new) element $e_{ab}$ that has no strict lower bound;
the strict upper bounds of $e_{ab}$ are just the elements of $\up^\sbX{\!\{a,b\}}$.  (Note:
$e_{ab}$ and $e_{ba}$ are the same element.)

Observe that if $\sbX$ has depth $n$, then so does $\widehat{\sbX}$, unless $n=2$ (in which
case $\widehat{\sbX}$ has depth $3$).  Also, since $X\in\textit{Up\/}(\widehat{\sbX})$,
we always have $\sbX^*\in\HHH({\widehat{\sbX}}^*)$.

The hat construction is illustrated below for
a poset $\sbP_6$ that will play a role in subsequent arguments.

{\small
\thicklines
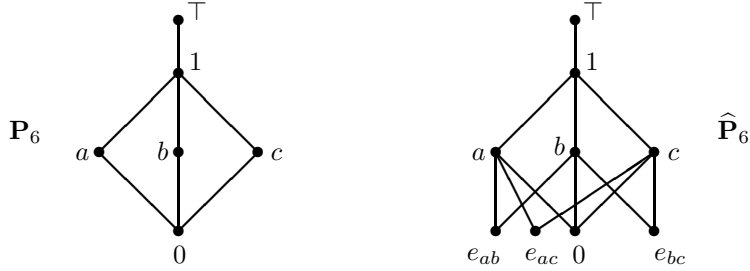
\begin{center}
\begin{picture}(0,110)(105,15)
\put(25,33){\circle*{4}}
\put(25,33){\line(1,1){30}}
\put(25,33){\line(0,1){80}}
\put(25,113){\circle*{4}}
\put(55,63){\circle*{4}}
\put(25,63){\circle*{4}}
\put(25,33){\line(-1,1){30}}
\put(-5,63){\circle*{4}}
\put(-5,63){\line(1,1){30}}
\put(25,93){\circle*{4}}
\put(55,63){\line(-1,1){30}}
\put(-39,68){$\mathbf{P}_6$}

\put(28,113){\small $\top$}
\put(29,94){\small $1$}
\put(60,60){\small $c$}
\put(-14,60){\small $a$}
\put(17,60){\small $b$}
\put(23,21){\small $0$}


\put(175,33){\circle*{4}}
\put(175,33){\line(1,1){30}}
\put(175,33){\line(0,1){80}}
\put(175,113){\circle*{4}}
\put(205,33){\circle*{4}}
\put(205,33){\line(0,1){30}}
\put(145,33){\line(1,1){30}}
\put(160,33){\line(-1,2){15}}
\put(160,33){\line(3,2){45}}
\put(205,33){\line(-1,1){30}}
\put(205,63){\circle*{4}}
\put(175,63){\circle*{4}}
\put(175,33){\line(-1,1){30}}
\put(145,33){\circle*{4}}
\put(145,33){\line(0,1){30}}
\put(160,33){\circle*{4}}
\put(145,63){\circle*{4}}
\put(145,63){\line(1,1){30}}
\put(175,93){\circle*{4}}
\put(205,63){\line(-1,1){30}}
\put(229,68){$\widehat{\mathbf{P}}_6$}

\put(178,113){\small $\top$}
\put(179,94){\small $1$}
\put(210,60){\small $c$}
\put(136,60){\small $a$}
\put(167,62){\small $b$}
\put(174,21){\small $0$}
\put(134,22){\small $e_{ab}$}
\put(156,22){\small $e_{ac}$}
\put(205,22){\small $e_{bc}$}

\end{picture}\nopagebreak



\end{center}
}


For each positive integer $n$, in the $n$-th direct power of the two-element chain,
let $\sbK_n$ be the subposet consisting of the least element, the greatest element, the $n$ atoms and
the $n$ co-atoms.   Note that $\sbK_n$ is bounded and has depth $3$.  Each atom of $\sbK_n$ is dominated by
just $n-1$ co-atoms, and each co-atom dominates just $n-1$ atoms.  Let
\[
\textup{$\mathcal{K}\seteq
$ the power set of
$\{\sbK_n: 3\leq n\in\mathbb{N}\}$,}
\]
so $\left|\mathcal{K}\right|=2^{\aleph_0}$.  Kuznetsov \cite{Kuz75}
proved that there are $2^{\aleph_0}$ distinct varieties of Brouwerian algebras of depth $3$,
by establishing the following:
\begin{equation}\label{kuznetsov}
\textup{for any $\class{C},\class{D}\in\mathcal{K}$, if $\class{C}\neq\class{D}$, then $\VVV(\class{C}^*)\neq\VVV(\class{D}^*)$,}
\end{equation}
where $\class{C}^*$ abbreviates $\{\sbX^*:\sbX\in\class{C}\}$.
(Ostensibly, \cite{Kuz75} deals with Heyting algebras, but its argument applies equally to
Brouwerian algebras.)


\begin{Lemma}\label{hard part}
Let\/ $\sbY,\sbZ,\sbW\in\{\sbP_6\}\cup\{\sbK_n:n\geq 3\}$\textup{,} where $\sbY\notin\{\sbP_6,\sbZ\}$\textup{.}
Then\/ $\sbY^*\notin\SSS\HHH({\widehat{\sbZ}}^*)$ and\/
${\sbW}^*\notin\III\SSS({\widehat{\sbZ}}^*)$\textup{.}
\end{Lemma}
\begin{proof}
%
%
%
%
Suppose $\sbY^*\in\SSS\HHH({\widehat{\sbZ}}^*)$.  Dualizing the injective/surjective homomorphisms,
we infer from Lemma \ref{1-1 and onto} that there exist $U\in\textit{Up\/}(\widehat{\sbZ})$ and a surjective
$p$-morphism $g\colon\sbU\mrig\sbY$.  Now $\sbY$ and $\widehat{\sbZ}$ have depth $3$, and $\sbY$
has a unique element of depth $3$, viz.\ its least element, $w$ say.
As $g$ is surjective, $w=g(u)$ for some $u\in U$.  Then $u$ has depth $3$, by (\ref{depth preservation}).
As $g$ is a $p$-morphism and $w$ has at least three distinct covers, each of depth $2$, the same is true
of $u$, by (\ref{p-morphism})
and (\ref{depth preservation}).  This prevents $u$ from having the form $e_{xy}$, so $u$ belongs to $\sbZ$.
As the least element of $\sbZ$ is its sole element of depth $3$, it is $u$.  Therefore, $Z\subseteq U$,
as $U\in\textit{Up\/}(\widehat{\sbZ})$.  Moreover, $\up{w}=\up{g(u)}=g\big[\!\up^{\widehat{\sbZ}}{u}\big]$,
i.e., $\sbY=g[\sbZ]$.

Then $\sbZ\neq\sbP_6$, because $\sbY$ has at least eight elements, while $\sbP_6$
has only six.  Thus, $\sbY,\sbZ$ are distinct elements of $\{\sbK_n:n\geq 3\}$.
As $g|_Z$ is a surjective $p$-morphism $\sbZ\mrig\sbY$, the homomorphism
$(g|_Z)^*\colon\sbY^*\mrig\sbZ^*$ is injective, by Lemma~\ref{1-1 and onto}, so $\sbY^*\in\III\SSS(\sbZ^*)$,
whence $\VVV(\sbY^*,\sbZ^*)=\VVV(\sbZ^*)$.  This contradicts (\ref{kuznetsov}),
because $\sbY\neq\sbZ$, so $\sbY^*\notin\SSS\HHH({\widehat{\sbZ}}^*)$.

Now suppose $\sbW^*\in\III\SSS({\widehat{\sbZ}}^*)$.  Then the situation in the first paragraph of the
present proof obtains, but with $\sbU=\widehat{\sbZ}$ and $\sbY=\sbW$.  Let $p,q,r$ be distinct covers of
$w$ in $\sbW$.  As we saw above, $u$ has distinct covers $p',q',r'$ (of depth $2$) that are mapped by
$g$ to $p,q,r$, respectively.
As $g$ is isotone, $g(e_{p'q'})$ is a common lower bound of the set $\{g(p'),g(q')\}=\{p,q\}$, so $g(e_{p'q'})=w$.
Then, because $w$ has three distinct covers (of depth $2$) in $\sbW$, it follows from (\ref{p-morphism}) and
(\ref{depth preservation}) that $e_{p'q'}$ has three distinct covers (of depth $2$) in $\widehat{\sbZ}$, but
this contradicts the definitions of $e_{p'q'}$ and $\widehat{\sbZ}$.  Thus, $\sbW^*\notin\III\SSS({\widehat{\sbZ}}^*)$.
\end{proof}

\begin{Lemma}\label{not sc}
Let\/ $\class{C}=\{{\widehat{\boldsymbol{P}}_{6}}\}\cup \{\widehat{\sbZ}:\sbZ\in\class{E}\}$\textup{,} where\/
$\class{E}\in\mathcal{K}$\textup{.}  Then the variety\/ $\VVV(\class{C}^*)$
is structurally incomplete.
\end{Lemma}
\begin{proof}
Let $\sbD$ be the direct product of the members of $\class{C}^*$.
Then $\VVV(\class{C}^*)=\VVV(\sbD)$, so it
suffices to show that $\VVV(\class{C}^*)\neq\QQQ(\sbD)$.
As
${\boldsymbol{P}_{6}}^{\ast} \in \HHH({\widehat{\boldsymbol{P}}_{6}}^*)\subseteq\VVV(\class{C}^*)$,
it is enough to prove that
${\boldsymbol{P}_{6}}^{\ast} \notin \QQQ(\sbD)$.
Observe that
\[
\textup{${\boldsymbol{P}_{6}}^{\ast} \in \QQQ(\sbD)$
 \ iff \
${\boldsymbol{P}_{6}}^{\ast} \in \III\SSS\PPU(\sbD)$
 \ iff \
${\boldsymbol{P}_{6}}^{\ast} \in \III\SSS(\sbD)$.}
\]
The first equivalence obtains because $\QQQ=\III\PPS\SSS\PPU$ and
${\boldsymbol{P}_{6}}^{\ast}$ is subdirectly irreducible; the second because
${\boldsymbol{P}_{6}}^{\ast}$ is finite and of finite type.
We must therefore show that
${\boldsymbol{P}_{6}}^{\ast}\notin\III\SSS(\sbD)$.

Suppose
${\boldsymbol{P}_{6}}^{\ast} \in \III\SSS(\sbD)$.
Then
${\boldsymbol{P}_{6}}^{\ast} \in \III\SSS\PPP(\class{C}^*)$.
As $\SSS\PPP=\PPS\SSS$, it follows (again from the subdirect irreducibility of
${\boldsymbol{P}_{6}}^{\ast}$)
that
${\boldsymbol{P}_{6}}^{\ast}$ embeds into ${\widehat{\sbZ}}^*$ for some $\sbZ\in\{\sbP_6\}\cup\class{E}$.
This contradicts Lemma~\ref{hard part}, so ${\boldsymbol{P}_{6}}^{\ast}\notin\III\SSS(\sbD)$.
\end{proof}

Every variety of Dunn monoids (e.g., Brouwerian algebras) has EDPC, by \cite[Thm.~3.55]{GJKO07}.
We shall require the following general result.

\begin{Theorem}\label{bp splitting}
\textup{(\cite{BP82},\,\cite[Thm.\,6.6]{Jon95})}
Let\/ $\class{K}$ be a variety of finite type, with EDPC, and let\/ $\sbA\in\class{K}$ be finite and
subdirectly irreducible.  Then there is a largest subvariety of\/ $\class{K}$ that excludes\/ $\sbA$\textup{.}
It consists of all\/ $\sbB\in\class{K}$ such that\/ $\sbA\notin\SSS\HHH(\sbB)$\textup{.}
\end{Theorem}

\begin{Lemma}\label{continuum lemma}
The set\/ $\{ \VVV(\class{C}^*):
\class{C}=\{{\widehat{\boldsymbol{P}}_{6}}\}\cup \{\widehat{\sbZ}:\sbZ\in\class{E}\} \text{ for some\/ }
\class{E}\in\mathcal{K}\}$
is a continuum of varieties of Brouwerian algebras of depth\/ $3$\textup{.}
\end{Lemma}

\begin{proof}
Suppose $\class{E}\cup\{\sbY\}\in\mathcal{K}$, where $\sbY\notin\{\sbP_6\}\cup\class{E}$.
It suffices to show that ${\widehat{\sbY}}^*\notin\VVV(\{{{\widehat{\sbP}}_6}^*\}\cup\{{\widehat{\sbZ}}^*:\sbZ\in\class{E}\})$.
As $\sbY^*\in\HHH({\widehat{\sbY}}^*)$,
it is enough to prove that
${\sbY}^*\notin\VVV(\{{{\widehat{\sbP}}_6}^*\}\cup\{{\widehat{\sbZ}}^*:\sbZ\in\class{E}\})$.

Let $\class{G}$ be the class of all Brouwerian algebras $\sbB$ such that
${{\sbY}}^*\notin \SSS\HHH(\sbB)$.  Since ${{\sbY}}^*$ is finite and subdirectly irreducible,
Theorem~\ref{bp splitting} shows that,
for any variety $\class{K}$ of
Brouwerian algebras, we have
${{\sbY}}^*\notin\class{K}$ iff $\class{K}\subseteq\class{G}$.
Therefore, it remains only to confirm that $\{{{\widehat{\sbP}}_6}^*\}\cup\{{\widehat{\sbZ}}^*:\sbZ\in\class{E}\}\subseteq\class{G}$,
i.e., that
${\sbY}^*\notin\SSS\HHH({{\widehat{\sbP}}_6}^*)$ and
${\sbY}^*\notin\SSS\HHH({{\widehat{\sbZ}}}^*)$ for all $\sbZ\in\class{E}$.
This is indeed the case, by Lemma~\ref{hard part}, because $\sbY\notin\{\sbP_6\}\cup\class{E}$.
\end{proof}

\begin{Theorem}\label{continuum brouwerian not sc}
The variety of Brouwerian algebras has\/ $2^{\aleph_0}$ structurally incomplete subvarieties (of depth $3$\textup{).}
\end{Theorem}
\begin{proof}
Use Lemmas~\ref{not sc} and \ref{continuum lemma}.
\end{proof}

This also completes the proof of Theorem~\ref{m continuum not sc}.
When we switch from Brouwerian to Heyting algebras, the duality theory undergoes slight definitional changes
(see \cite[Sec.\,4]{BMR17}, for instance).
The proof of Theorem~\ref{continuum brouwerian not sc} adapts easily, however, and the result remains true
when we replace `Brouwerian' by `Heyting' in its statement.

In a quasivariety $\class{K}$, a quasi-equation
\begin{equation*}
\left(\varphi_1\thickapprox\psi_1\aand\dots\aand\varphi_n\thickapprox\psi_n\right)\Longrightarrow\varphi\thickapprox\psi
\end{equation*}
is said to be \emph{admissible} provided that, for every substitution $h$,
\[
\textup{if \,$\class{K}\models h(\varphi_i)\thickapprox h(\psi_i)$\, for $i=1,\dots,n$, \,then
\,$\class{K}\models h(\varphi)\thickapprox h(\psi)$.}
\]
By Theorem~\ref{free}, $\class{K}$ is SC iff it satisfies its own admissible quasi-equations.

Mints \cite{Min76} showed (in effect) that the variety $\class{BRA}$ of \emph{all} Brouwerian algebras
is not SC, by proving that the following quasi-equation (not satisfied by $\class{BRA}$) is admissible in $\class{BRA}$:
\begin{equation}\label{mints rule}
x\rig y\leqslant x\vee z \,\Longrightarrow\, ((x\rig y)\rig
x)\,\vee\,((x\rig y)\rig z)\thickapprox e.
\end{equation}

For each term $\varphi$ over $\class{BRA}$,
we define a term $\varphi^\diamond$ in the following recursive
manner: $e^\diamond=e$\,; $x^\diamond=x\wedge e$ ($x\in \mathit{Var}$); $(\al*\be)^\diamond=\al^\diamond * \be^\diamond$
($*\in\{\wedge,\vee\}$); $(\al\rig\be)^\diamond=(\al^\diamond\rig\be^\diamond)\wedge e$.  In this signature,
the \emph{amendment} $\Phi^\diamond$ of a quasi-equation $\Phi$ results from replacing each term $\varphi$ occurring in $\Phi$
with the term $\varphi^\diamond$.
Using ideas of Iemhoff \cite{Iem01}, one can prove that if a quasi-equation $\Phi$ is admissible in $\class{BRA}$, then $\Phi^\diamond$ is admissible in the
variety $\class{DM}$ of all Dunn monoids; moreover, if $\class{BRA}\not\models\Phi$, then $\class{DM}\not\models\Phi^\diamond$.
This implies, of course, that $\class{DM}$ is not SC.

It is tempting to try to extend the argument from $\class{DM}$
to the variety $\class{M}$ of Definition~\ref{m defn}, using reflections.  Unfortunately, however, it turns out that
the amendment of (\ref{mints rule}) is not admissible in $\class{M}$.

\medskip

{\small
\noindent {\bf Acknowledgment.}  The first author thanks Miguel Campercholi for stimulating conversations on this topic.}

\end{document}